\documentclass{amsart}
\usepackage{color}
\usepackage{graphicx}
\newtheorem{theorem}{Theorem}[]
\definecolor{brightred}{rgb}{1.0,0.9,0.9}
\definecolor{brightblue}{rgb}{0.9,0.9,1.0}
\definecolor{brightyellow}{rgb}{1.0,1.0,0.8}
\definecolor{brightgreen}{rgb}{0.9,1.0,0.9}
\def\satz#1#2{\begin{center}\fcolorbox{brightyellow}{brightyellow}{\parbox{9cm}{\begin{theorem}[#1]#2 \end{theorem}}}\end{center}}
\def\df#1{{\bf \underline{#1}}}
\title{The theorems of Green-Stokes,Gauss-Bonnet and Poincar\'e-Hopf in Graph Theory}
\date{January 27, 2012}
\author{Oliver Knill}
\email{knill@math.harvard.edu}
\address{
        Department of Mathematics \\
        Harvard University \\
        Cambridge, MA, 02138
        }

\begin{document}
\maketitle

\begin{abstract}
By proving graph theoretical versions of 
Green-Stokes, Gauss-Bonnet and Poincar\'e-Hopf, core ideas of undergraduate mathematics 
can be illustrated in a simple graph theoretical setting. In this pedagogical exposition
we present the main proofs on a single page and add illustrations. 
While discrete Stokes is at least 100 years old, the other two results for graphs 
were found only recently. 
\end{abstract}

\section{Definitions}

A \df{simple graph} $G=(V,E)$ is a finite set $V$ paired with a finite set $E$ of unordered pairs $e=\{ a,b \}$ with $a,b \in V$. 
A simple graph has no multiple connections and no self-loops: every $e \in E$ appears only once and no $\{ a,a \}$ is in $E$. 
Elements in $V$ are called \df{vertices}, elements in $E$ are called \df{edges}. 
Given a simple graph $G=(V,E)$, denote by $S(v)$ the \df{unit sphere} of a vertex $v$. It is a subgraph generated by the set of vertices
directly connected to $v$. Denote by $G_k$ the set of complete $K_{k+1}$ subgraphs of $G$. Elements in $G_k$ are also called \df{cliques}. 
The set $G_2$ for example is the set of all triangles in $G$. Of course, $G_0=V$ and $G_1=E$. If the cardinality of $G_k$ is 
denoted by $v_k$, the \df{Euler characteristic} of $G$ is defined as
$\chi(G) = \sum_{k=0}^{\infty} (-1)^k v_k$, a finite sum. For example, if no tetrahedral subgraphs $K_4$ exist
in $G$, then $\chi(G)=v-e+f$, where $v=v_0={\rm ord}(G)$ is the \df{order}, the number of vertices, $e=|G|=v_1$ is the \df{size}, 
the number of edges and $f=v_2$ the number of \df{triangles} $K_3$. 
\df{Dimension} is defined inductively as ${\rm dim}(G) = 1+|V|^{-1} \sum_{v \in V} {\rm dim}(S(v))$ with 
${\rm dim}(\emptyset)=-1$. Cyclic graphs, trees or the dodecahedron are examples of graphs of dimension $1$, 
a triangle $K_3$, an octahedron or icosahedron has dimension $2$. A tetrahedron has dimension $3$. 
A complete graph $K_{k+1}$ on $k+1$ vertices has dimension $k$. 
Dimension is defined for any graph but can become a fraction. For a truncated cube $G$ for example, each unit
sphere $S(v)$ is a graph of $3$ vertices and one edge, a graph of dimension ${\rm dim}(S(v))=2/3$ so that ${\rm dim}(G)=5/3$. 
The Euler characteristic of this graph $G$ is $\chi(G) = v-e+f = 24-36+8=-4$.
A \df{$k$-form} is a function on $G_k$ which is antisymmetric in its $(k+1)$ arguments. The set $\Omega^k$ of all $k$-forms 
is a vector space of dimension $v_k$. The remaining sign ambiguity can be fixed by introducing an orientation on the 
graph: a $K_{k+1}$ subgraph is called a \df{maximal simplex} if it is not contained in a larger $K_l$ graph.
An \df{orientation} attaches a $k$-form $m$ to each maximal simplex with value $1$. It induces forms on smaller dimensional faces. 
If $m$ cancels on intersections of maximal graphs, it is a "volume form" and $G$ is called \df{orientable}. An icosahedron for 
example has triangles as maximal simplices.
It is orientable. A wheel graph $W_6$ in which two opposite edges are identified models a M\"obius strip and is not orientable.
A $0$-form is a function on $V=G_0$ and also called a \df{scalar function}. 
Call $df(a,b)=f(b)-f(a)$ the \df{gradient}. It is defined as a $1$-form if $G$ has an orientation. 
Without an orientation, we can still look at the \df{directional derivative} $D_ef(a) = f(b)-f(a)$ if $e=\{a,b\} \in E$ 
is an edge attached to $a$. Define the \df{exit set} $S_f^-(v)=\{ w \in S(v) \; | \; f(w)-f(v)<0  \; \}$
and the \df{index} $i_f(v) = 1-\chi(S_f^-(v))$. A vertex $v$ is a \df{critical point} if $i_f(v) \neq 0$.
If $k \geq 1$ and $G$ is $k$-dimensional, a vertex $v$ is an \df{interior point} if $S(v)$ is a $(k-1)$-dimensional 
graph for which every point is an interior point within $S(v)$; for $k \geq 2$ we ask $S(v)$ to be connected. 
The base induction assumption is that an interior point of a one-dimensional graph has two neighbors.
A vertex $v$ of a $k$-dimensional graph $G$ is a \df{boundary point} if 
$S(v)$ is a $(k-1)$-dimensional (for $k \geq 2$  connected) graph in which every vertex is either a 
boundary or interior point and both are not empty.  The seed assumption is that for $k=1$, the graph $S(v)$ has one vertex.
A $k$-dimensional graph $H$ is a \df{graph with boundary} $dG$ if
every $v \in V$ is an interior point or a boundary point. Glue two copies of $H$ along the boundary gives a graph $G$ without boundary.
A wheel graph $W_k$ is an example of a $2$-dimensional 
graph with boundary if $k \geq 4$. The boundary is the cyclic one dimensional graph $C_k$. Cut an octahedron in two gives $W_4$.
For an oriented graph $G$, the \df{exterior derivative} $d: \Omega_{k} \to \Omega_{k+1}$ is defined as
$df(x)=\sum_i (-1)^i f(x_0,\dots,\hat{x_i},\dots,x_k)$, where $\hat{x}$ denotes a variable taken away. 
For example $df(x,y,z) = f(y,z) - f(x,z) + f(x,y)$ is a function on triangles 
called the \df{curl} of a $1$-form $f$. A form is \df{closed} if $df=0$. It is \df{exact} if $f=dg$. The vector space
$H^k(G)$ of closed forms modulo exact forms is a \df{cohomology group} of dimension $b_k$, the \df{Betti number}. Example: $b_0$ is
the number of \df{connected components}. The \df{cohomological Euler characteristic} is $\sum_{k=0}^{\infty} (-1)^k b_k$. 
For a $k$-form define the \df{integral} $\sum_G f = \sum_{v \in G_k} f(v)$. 
Let $V_k(v)$ be the number of $K_{k+1}$ subgraphs of $S_k(v)$. Especially, $V_0(v)$ is the \df{degree} ${\rm deg}(v)$ of $v$, the
order of $S(v)$. The local quantity $K(v) = \sum_{k=0}^{\infty} (-1)^k V_{k-1}(v)/(k+1)$ is 
called the \df{curvature} of the graph at $v$.  The sum is of course finite. For a $2$-dimensional graph without boundary, where
$S(v)$ has the same order and size, it is $K(v) = 1-{\rm ord}(S(v))/2+|S(v)|/3 = 1-|S(v)|/6=1-{\rm deg}(v)/6$.
For a 1-dimensional graph with or without boundary and trees in particular, $K(v) = 1-{\rm ord}(S(v))/2=1-{\rm deg}(v)/2$. 

\section{Theorems}

For an arbitrary finite simple graph we have \cite{cherngaussbonnet}

\satz{Gauss-Bonnet}{
$\sum_{v \in V} K(v) = \chi(G)$.
\label{theorem1}
}

For an arbitrary finite simple graph and injective $f:V \to \mathbf{R}$, we have \cite{poincarehopf}

\satz{Poincare-Hopf}{
$\sum_{v \in V} i_f(v) = \chi(G)$
\label{theorem2}
}

Assume $f$ is a $(k-1)$-form and $G$ is an oriented $k$-dimensional graph with boundary, then
\footnote{This is combinatorial topology known at the time of Poincar\'e. 
Graph theoretical versions have since appeared in many dialects in numerical, physical or computer science contexts:
take any theory for simplicial complexes and remove the Euclidean "fillings". Samples are \cite{GSD,Regge,MarsdenDesbrun}.
Discrete differential geometry comes in many flavors: 
examples are computational geometry \cite{Devadoss,ComputationalGeometry}, integrable systems 
\cite{BobenkoSuris}, computer graphics or computational-numerical methods (i.e. \cite{Bobenko,CompElectro2002}.
New here for Stokes is only the notion of {\bf dimension} \cite{elemente11}
which allows to formulate the result conveniently from
{\bf within graph theory}. For graphs $G$ with boundary $dG$, the later remains a graph. In general it is 
only a \df{chain}, an element in the group of integer valued functions on ${\mathcal{G}}=\bigcup_k G_k$ 
usually written as $\sum_{c \in {\mathcal{G}}} a_c c$. }

\satz{Green-Stokes}{
$\sum_{G} df = \sum_{dG} f$.
\label{theorem3}
}


\section{Proofs}

\fcolorbox{brightblue}{brightblue}{Proof of \ref{theorem1}:} the {\bf transfer equations} are
\fcolorbox{brightgreen}{brightgreen}{$\sum_{v \in V} V_{k-1}(v) = (k+1) v_k$}.
By definition of curvature, we have 
$$  \sum_{v \in V} K(v) = \sum_{v \in V} \sum_{k=0}^{\infty} (-1)^k \frac{V_{k-1}(v)}{k+1} \; . $$
Since the sums are finite, we can change the order of summation. Using the transfer equations we get
$$  \sum_{v \in V} K(v) = \sum_{k=0}^{\infty} \sum_{v \in V} (-1)^k \frac{V_{k-1}(v)}{k+1}  
                        = \sum_{k=0}^{\infty} (-1)^k v_k = \chi(G) \;. $$

\fcolorbox{brightblue}{brightblue}{Proof of \ref{theorem2}:}
the number of $k$ simplices $V_k^-(v)$ in the exit set $S^-(v)$ 
and the number of $k$ simplices $V_k^+(v)$ in the entrance set $S^+(v)$ are complemented within $S(v)$ 
by the number $W_k(v)$ of $k$ simplices which contain both vertices from $S^-(v)$ and $S^+(v)$. By definition,
$V_k(v) = W_k(v) + V_k^+(v) + V_k^-(v)$. 
The index $i_f(v)$ is the same for all injective functions $f:V \to \mathbf{R}$. 
The {\bf intermediate equations}
are \fcolorbox{brightgreen}{brightgreen}{$\sum_{v \in V} W_k(v) = k v_{k+1}$}.
Let  $\chi'(G) = \sum_v i_f(v)$. Because replacing $f$ and $-f$ switches $S^+$ with $S^-$ and the sum is the same, we
can prove $2 v_0 - \sum_{v \in V} \chi(S^+(v)) + \chi(S^-(v)) = 2 \chi'(G)$ instead. 
The transfer equations and intermediate equations give
\begin{eqnarray*}
 2 \chi'(G) &=& 2v_0 + \sum_{k=0}^{\infty} (-1)^k \sum_{v \in V} (V_k^-(v) + V_k^+(v)) 
             =  2v_0 + \sum_{k=0}^{\infty} (-1)^k \sum_{v \in V} (V_k(v)   - W_k(v)  ) \\
            &=& 2v_0 + \sum_{k=0}^{\infty} (-1)^k [ (k+2) v_{k+1} - k v_{k+1}]  
             =  2v_0 + \sum_{k=0}^{\infty} (-1)^k 2 v_{k+1} = 2 \chi(G) \; .
\end{eqnarray*}

\fcolorbox{brightblue}{brightblue}{Proof of \ref{theorem3}:} denote a $k$-simplex graph $K_k$ by $(x_0,...,x_k)$. From
$df(x_0,...,x_n) = \sum_k (-1)^k f(x_0,...,\hat{x}_k,...,x_n)$ and algebraic boundary
$d[x_0,...,x_n] = \sum_k (-1)^k (x_0, ...,\hat{x}_k,...,x_n))$,
Stokes theorem is obvious for a single simplex:
$$ \sum_{G} df = \sum_k (-1)^k f(x_0, \dots ,\hat{x}_k, \dots ,x_n) = \sum_{dG} f  \; . $$
Gluing $k$-dimensional simplices cancels boundary.
A $k$-dimensional graph with boundary is a union of $k$-dimensional simplices identified along $(k-1)$-
dimensional simplices. A $k$-dimensional oriented graph with boundary can be built by gluing cliques
as long as the orientation $k$-form can be extended. We also used that the {\bf boundary as a graph} agrees
with the {\bf algebraic boundary} if differently oriented boundary pieces cancel.

\pagebreak

\section{Illustrations} 

\parbox{16.8cm}{
\parbox{6cm}{\scalebox{0.28}{\includegraphics{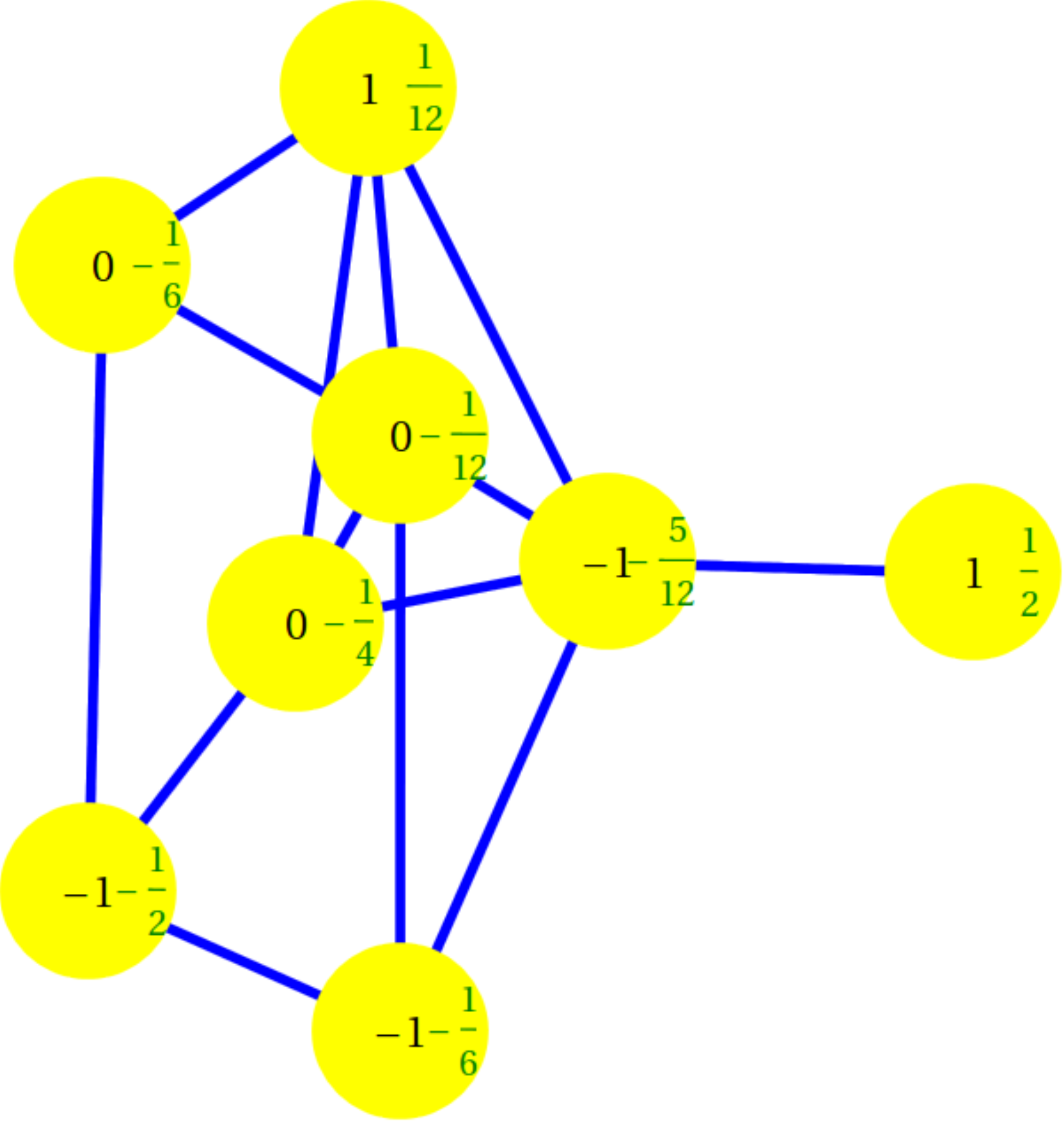}} }
\parbox{10cm}{
{\bf Illustration to Theorem 1.}
The figure shows a graph $G=(V,E)$ with $v_0=|V|=8$ vertices and $v_1=|E|=14$ edges. The Euler characteristic is
$\chi(G) = -1$. To describe the unit spheres, we use the notation $G \cup H$ for
a union of two disconnected graphs, $P_k$ is the $0$ dimensional graph with $k$ vertices
and no edges, $K_k$ is the complete graph with $k$ edges and $K_k^+$ is a complete graph with
an additional vertex and an edge connecting the appendix to $K_k$. 
The curvatures $K(v)$ are indicated near each vertex. We only compute of half of them.
If we add the curvatures up we get $-1$. 
\begin{center}
\begin{tabular}{|l|lllll|} \hline
  $v$ & $S(v)$          & $V_0(v)$       & $V_1(v)$  & $V_2(v)$ & $K(v)$             \\ \hline \hline
  $1$ & $P_1$           &  1             & 0         & 0        & 1-1/2+0/3-0/4=1/2  \\ \hline
  $5$ & $P_1\cup K_3^+$ &  5             & 4         & 1        & 1-5/2+4/3-1/4=-5/12\\ \hline
  $7$ & $P_3$           &  3             & 0         & 0        & 1-3/2+0/3+0/4=-1/2 \\ \hline
  $8$ & $P_1\cup K_2$   &  3             & 1         & 0        & 1-3/2+1/3+0/4=-1/6  \\ \hline
\end{tabular}
\end{center}
}
}

\vspace{1cm}

\parbox{16.8cm}{
\parbox{6cm}{\scalebox{0.28}{\includegraphics{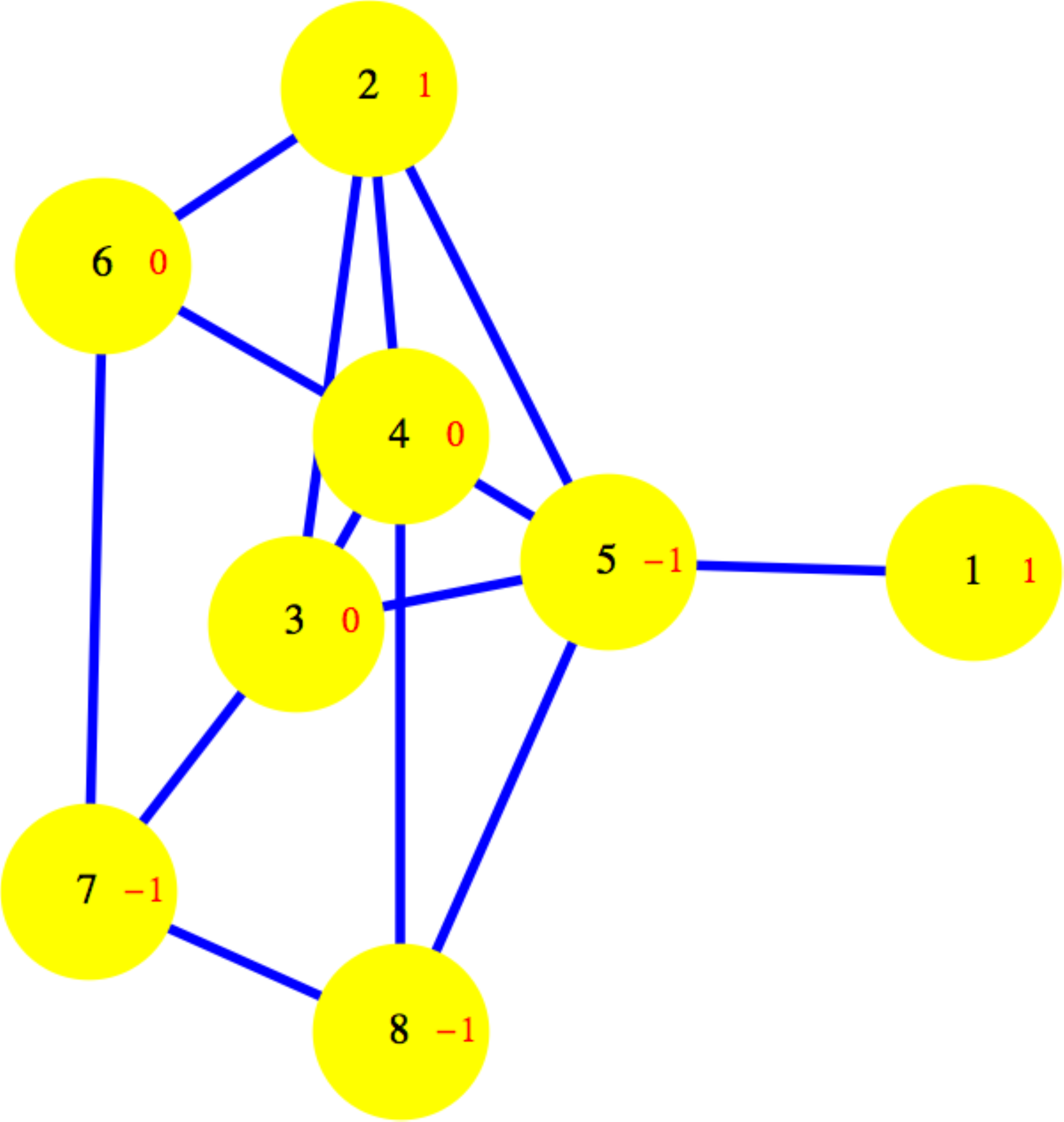}} }
\parbox{10cm}{
{\bf Illustration to Theorem 2.}
A graph $G=(V,E)$ with $v_0=|V|=8$ vertices and $v_1=|E|=14$ edges. 
There are $v_2=6$ triangles and $v_4=1$ tetrahedra. 
The Euler characteristic is $\chi(G) = 8-14+6-1=-1$. The graph in the picture is
equipped also with a Morse function which is used to label the vertices.
The indices are noted near each vertex. Lets compute the indices $i_f(v)$:
\begin{center}
\begin{tabular}{|l|llll|} \hline
  $v$ & $S(v)$        & $S^-(v)$       & $\chi(S^-(v))$  & $i(v)$   \\ \hline \hline
  $1$ & $P_1$         & $\emptyset$    & $0$ & $1$ \\ \hline
  $5$ & $P_1\cup K_3^+$ & $P_1 \cup K_3$ & $2$ & $-1$ \\ \hline
  $7$ & $P_3$         & $P_2$          & $2$ & $-1$ \\ \hline
  $8$ & $P_1\cup K_2$ & $P_1 \cup I_4$ & $2$ & $-1$ \\ \hline
\end{tabular}
\end{center}
}
}

\vspace{1cm}

\parbox{16.8cm}{
\parbox{6cm}{\scalebox{0.28}{\includegraphics{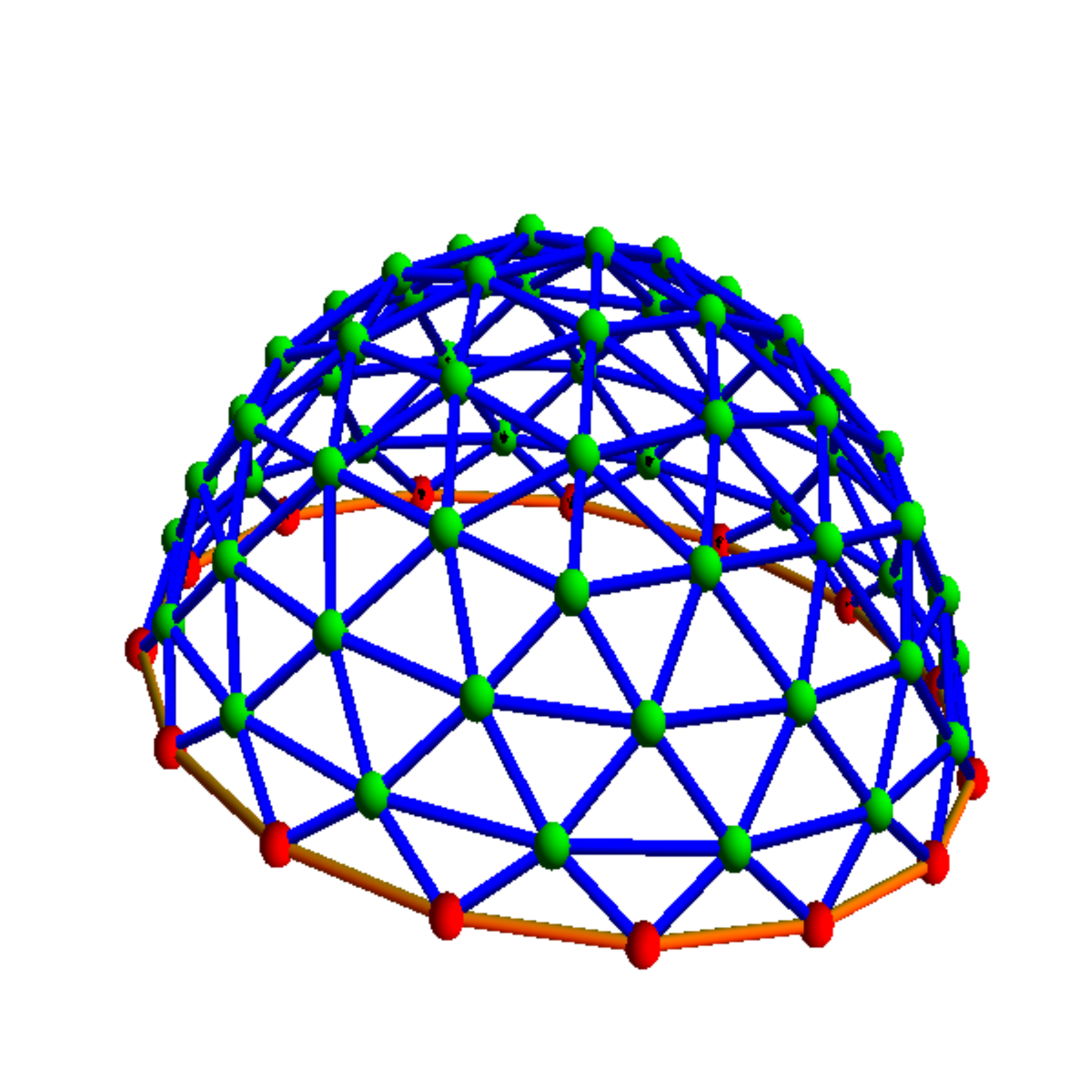}} }
\parbox{10cm}{
{\bf Illustration to Theorem 3.}
The figure shows a geodesic dome $G$. It is a two dimensional graph with boundary. The
graph has $v=76$ vertices and $e=210$ edges and $f=135$ triangles and
$\chi(G) = v-e+f=1$.  There are 61 interior points, where the unit sphere is a circle and 
15 boundary points for which the unit sphere is an interval graph. 
The dome has 6 vertices with degree $5$ and so curvature $1/6$. The sum of curvatures is $1$.
The graph is orientable.  The maximal simplices are the triangles. An orientation assigns a 
permutation on each triangle. We can assume for example that
each triangle is oriented in a counter clockwise way when looking from above onto
the dome. A $1$-form $f$ is a function on edges. The {\bf curl} $df$ is a two
form. It is a function on the triangles. The value of $df$ is obtained by 
adding up the three values of $f$ along the boundary of the triangle. 
The sum of the curls cancel in the interior. Only the boundary part 
survives. The sum $\sum_{v \in dG} f$ is the "line integral" of $f$ along 
a circular graph $C_{15}$. You see Stokes theorem. 
}
}

\pagebreak

Here are families of graphs, where the curvature is indicated at every vertex: 

\parbox{16.8cm}{
\parbox{5.5cm}{\scalebox{0.25}{\includegraphics{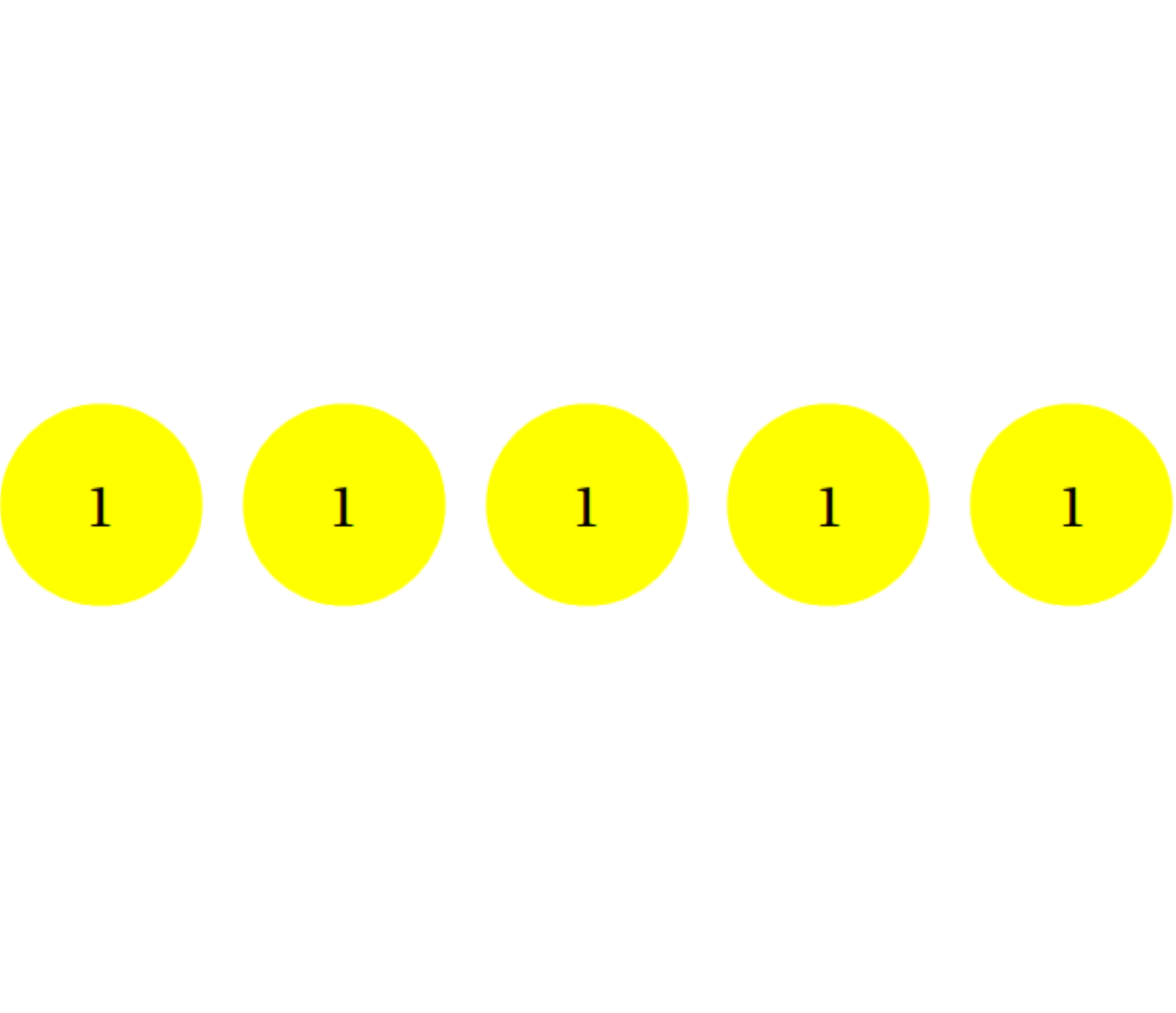}}}
\parbox{5.5cm}{\scalebox{0.25}{\includegraphics{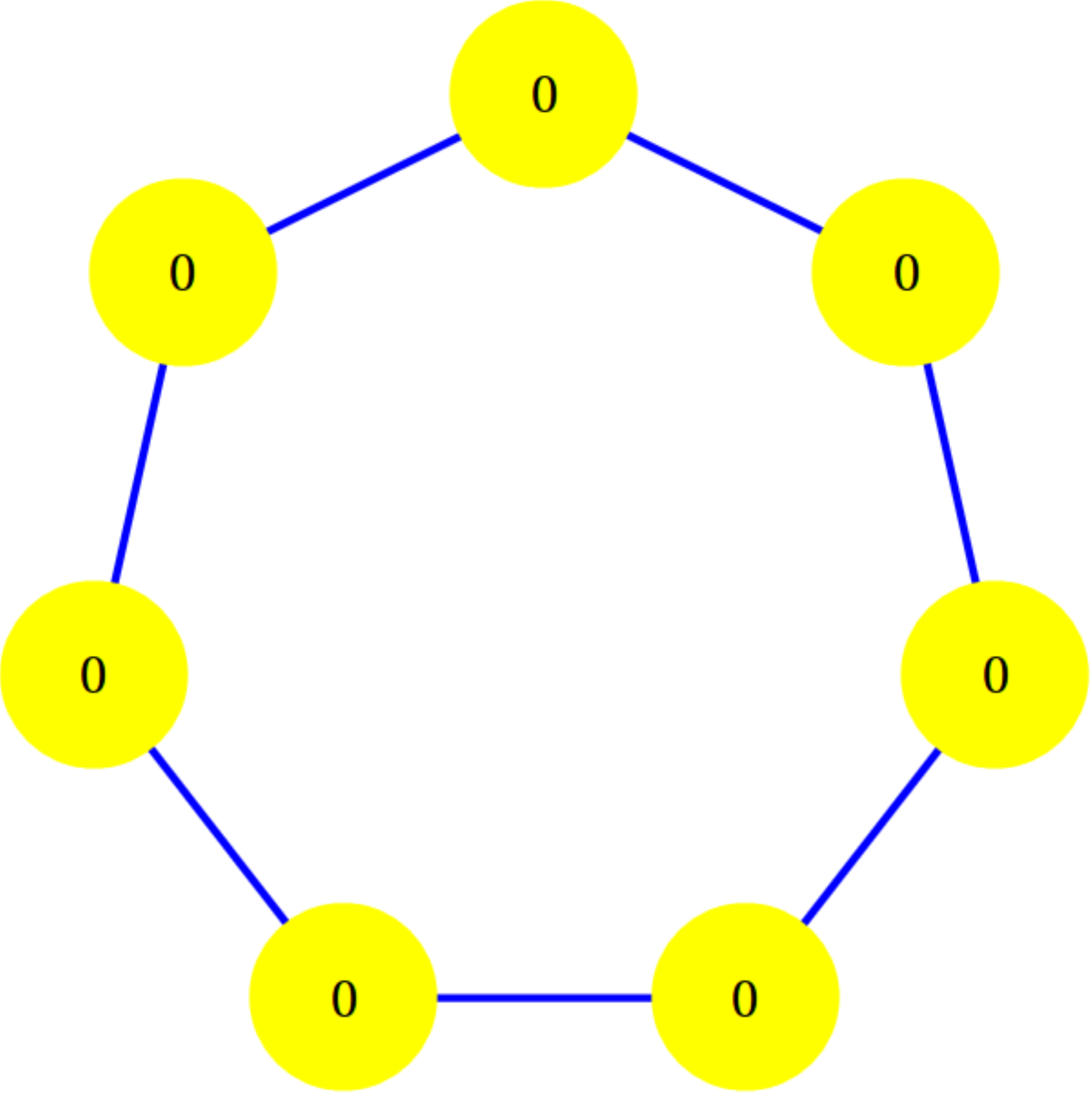}}}
\parbox{5.5cm}{\scalebox{0.25}{\includegraphics{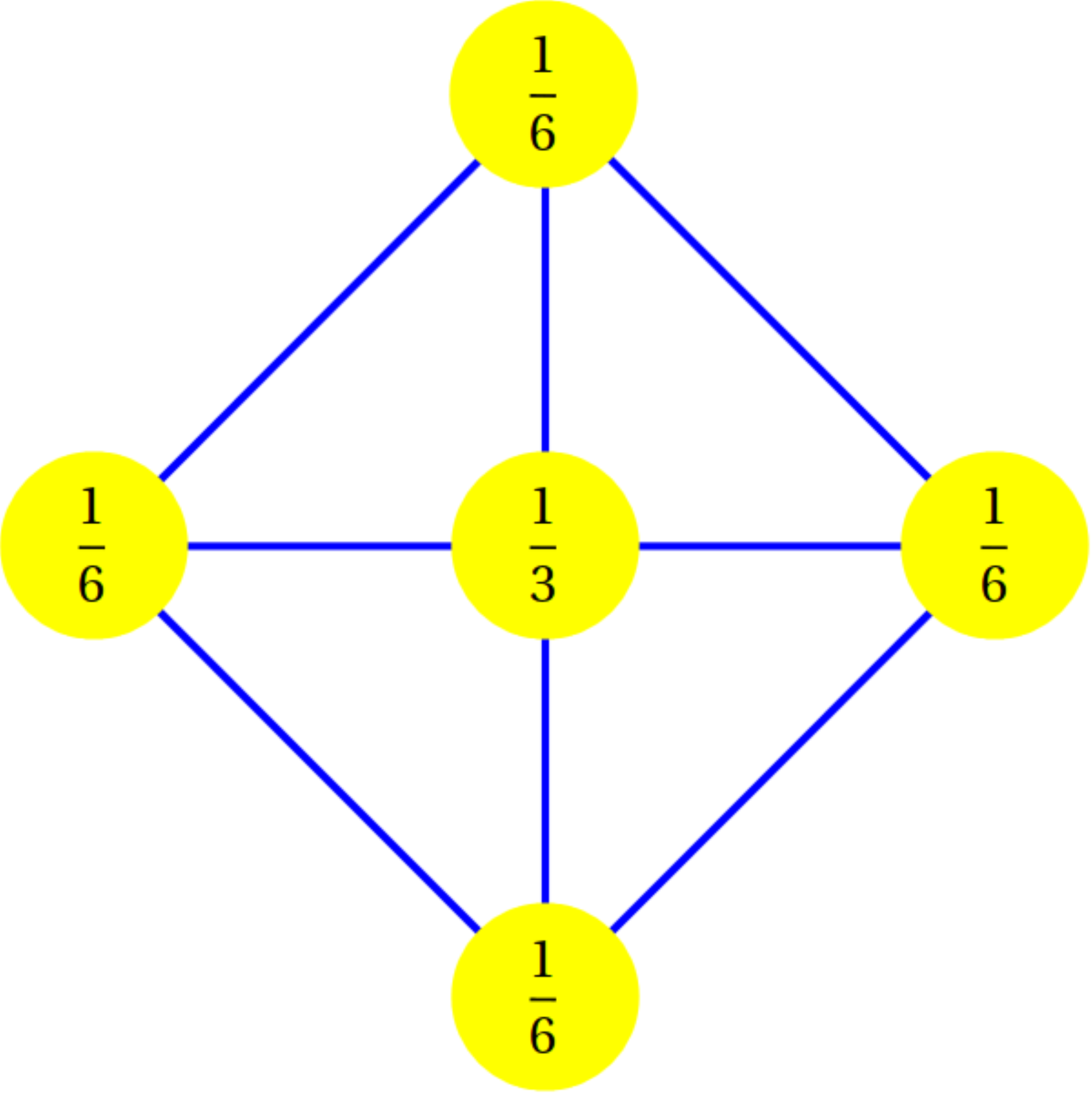}}}
}

\parbox{16.8cm}{ \parbox{5.5cm}{ $P_5$ } \parbox{5.5cm}{ $C_7$ } \parbox{5.5cm}{$W_4$ } }

\parbox{16.8cm}{
\parbox{5.5cm}{\scalebox{0.25}{\includegraphics{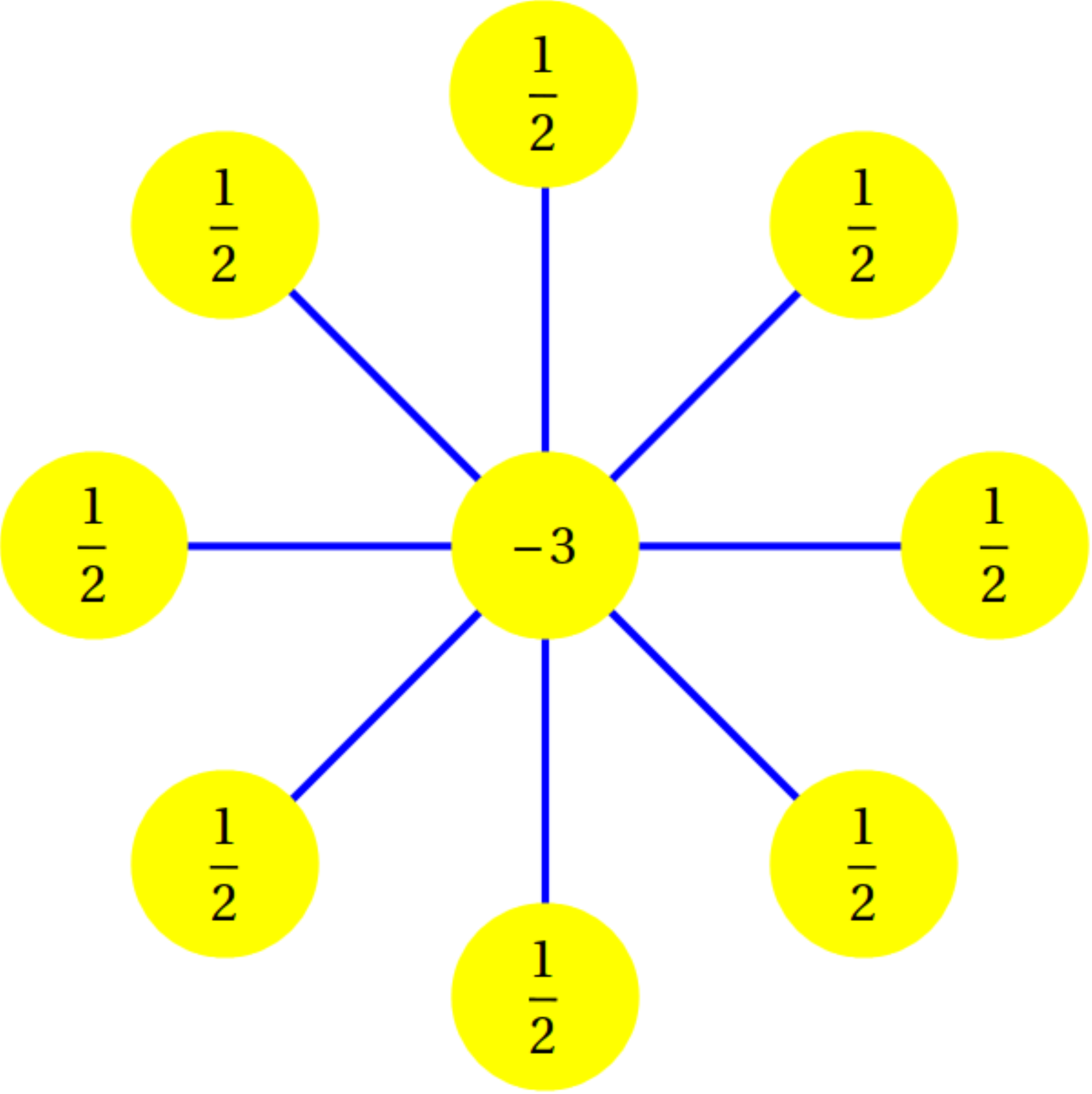}}}
\parbox{5.5cm}{\scalebox{0.25}{\includegraphics{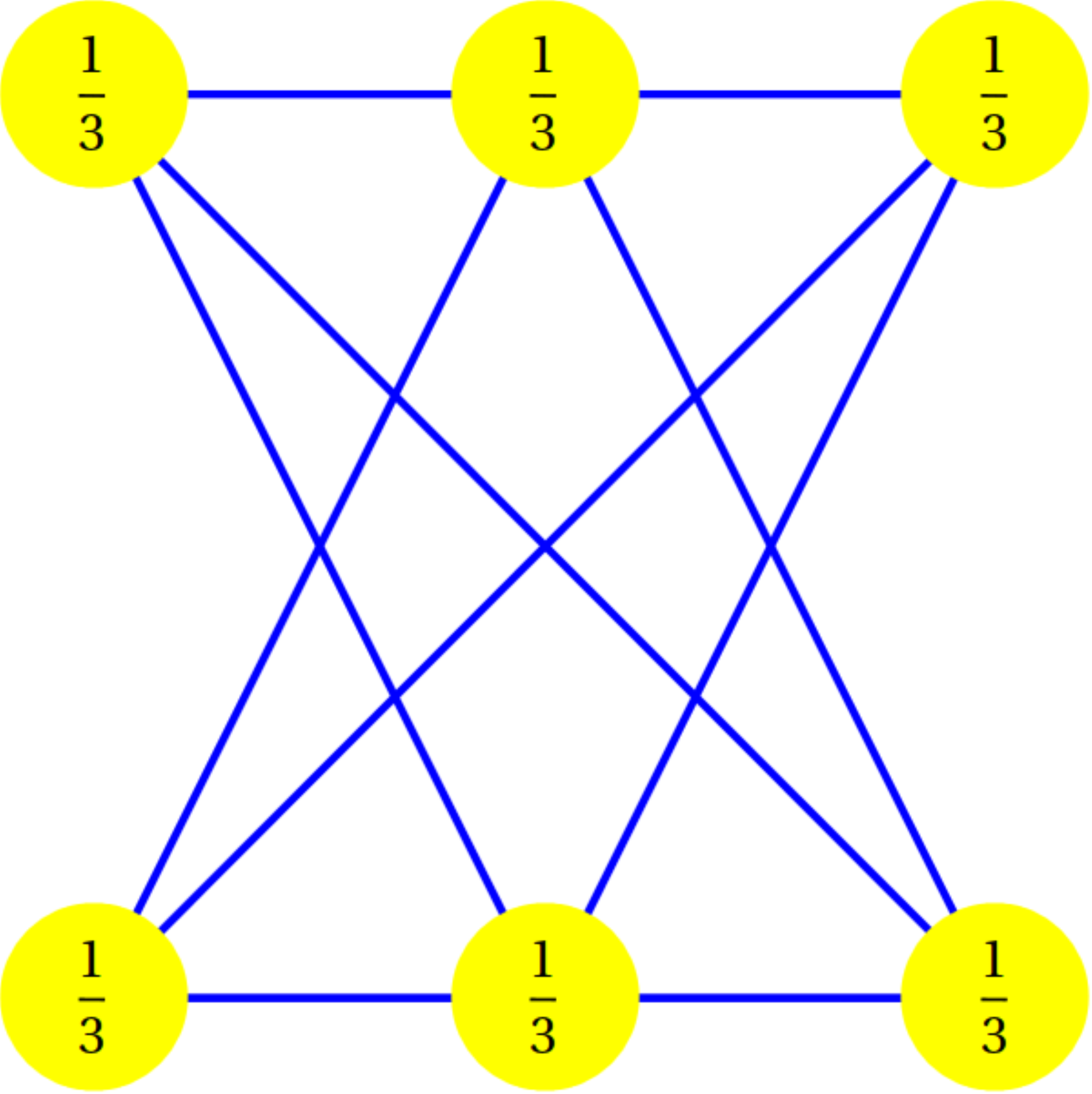}}}
\parbox{5.5cm}{\scalebox{0.25}{\includegraphics{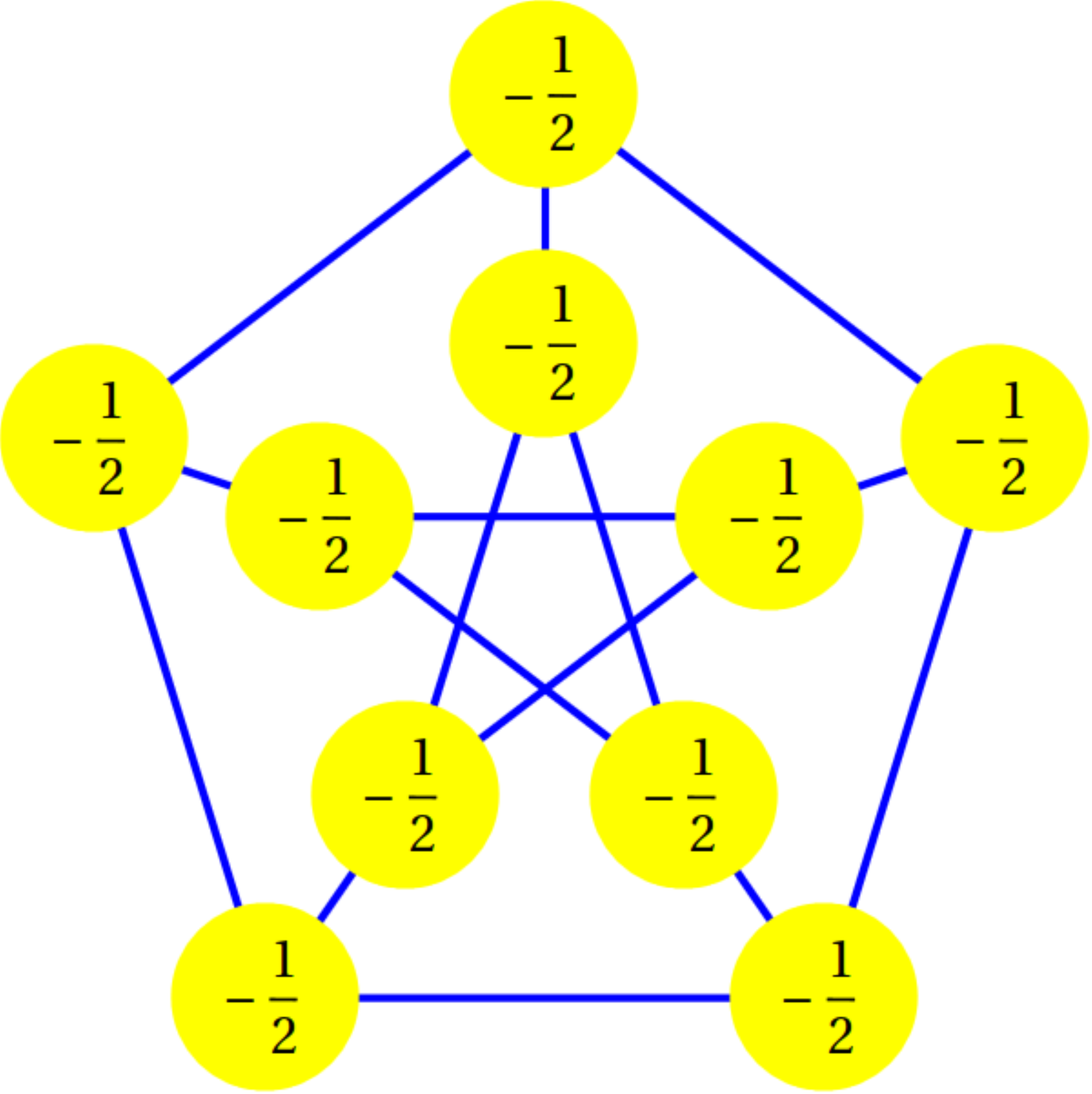}}}
}

\parbox{16.8cm}{ \parbox{5.5cm}{ $S_8$ } \parbox{5.5cm}{ $T_6$ } \parbox{5.5cm}{ $P_{5,2}$ } }

\parbox{16.8cm}{
\parbox{5.5cm}{\scalebox{0.25}{\includegraphics{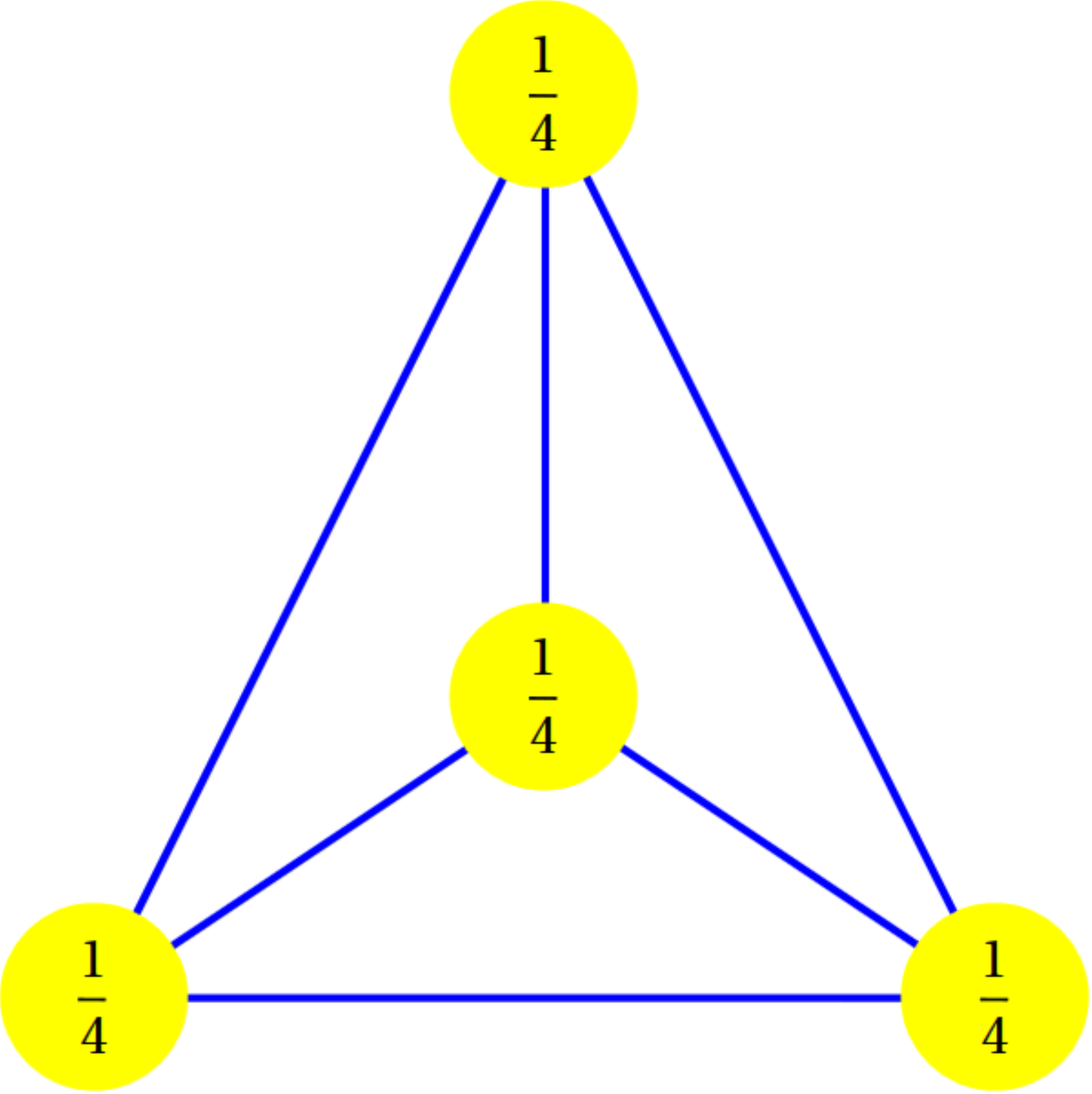}}}
\parbox{5.5cm}{\scalebox{0.25}{\includegraphics{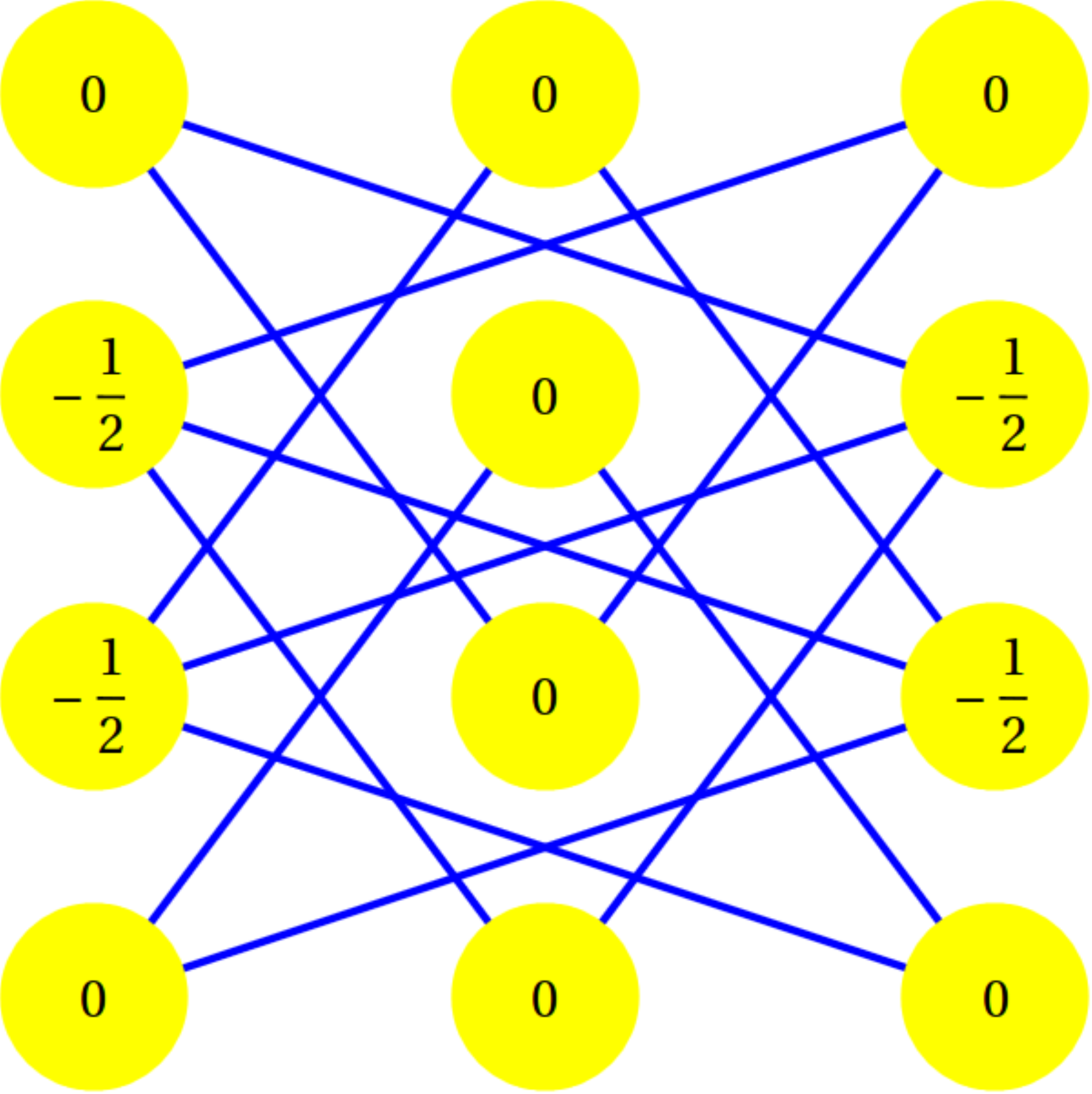}}}
\parbox{5.5cm}{\scalebox{0.25}{\includegraphics{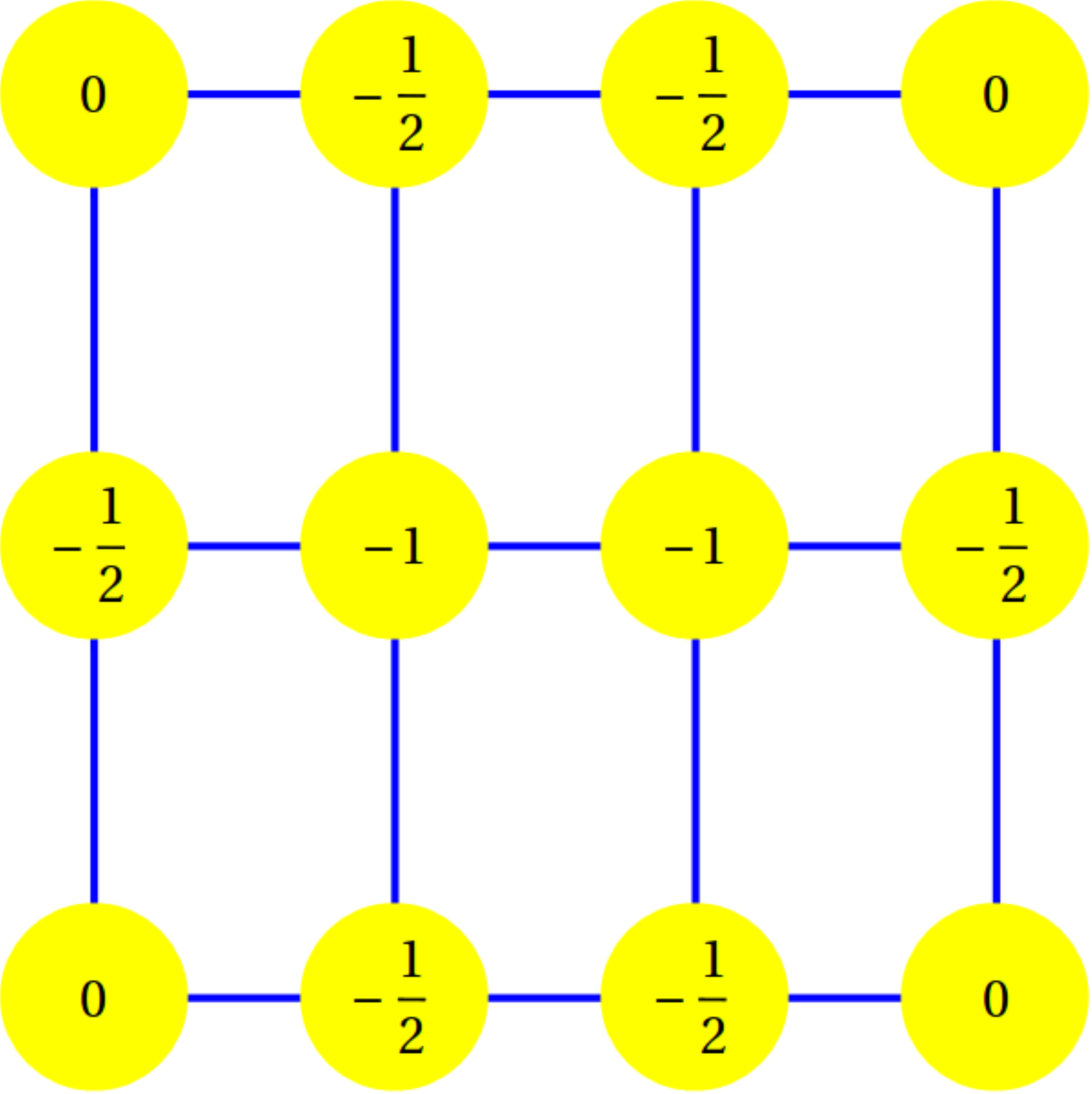}}}
}

\parbox{16.8cm}{ \parbox{5.5cm}{ $K_4$ } \parbox{5.5cm}{ $K_{3,4}$ } \parbox{5.5cm}{ $G_{3,4}$ } }

\parbox{16.8cm}{
\parbox{5.5cm}{\scalebox{0.25}{\includegraphics{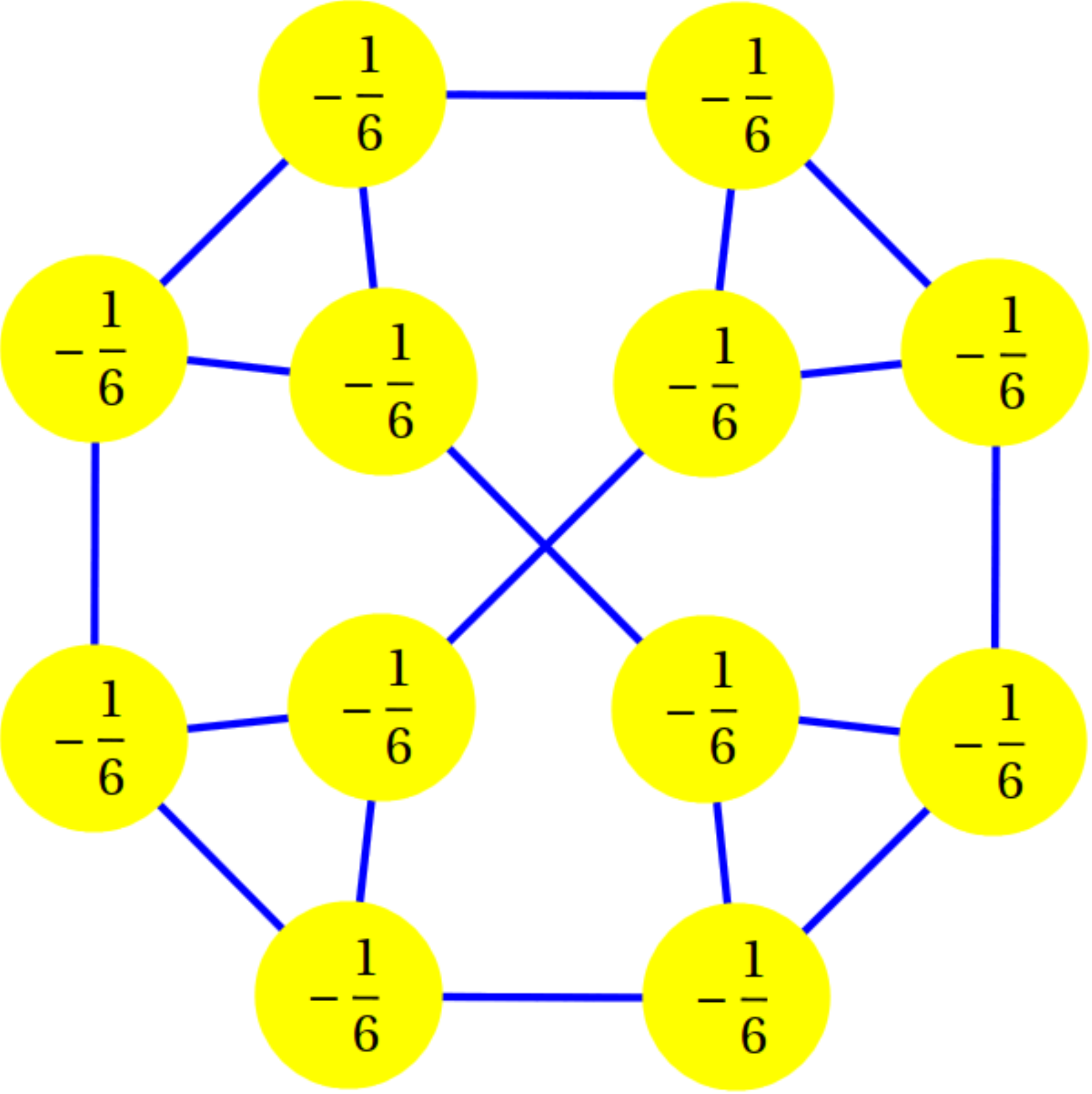}}}
\parbox{5.5cm}{\scalebox{0.25}{\includegraphics{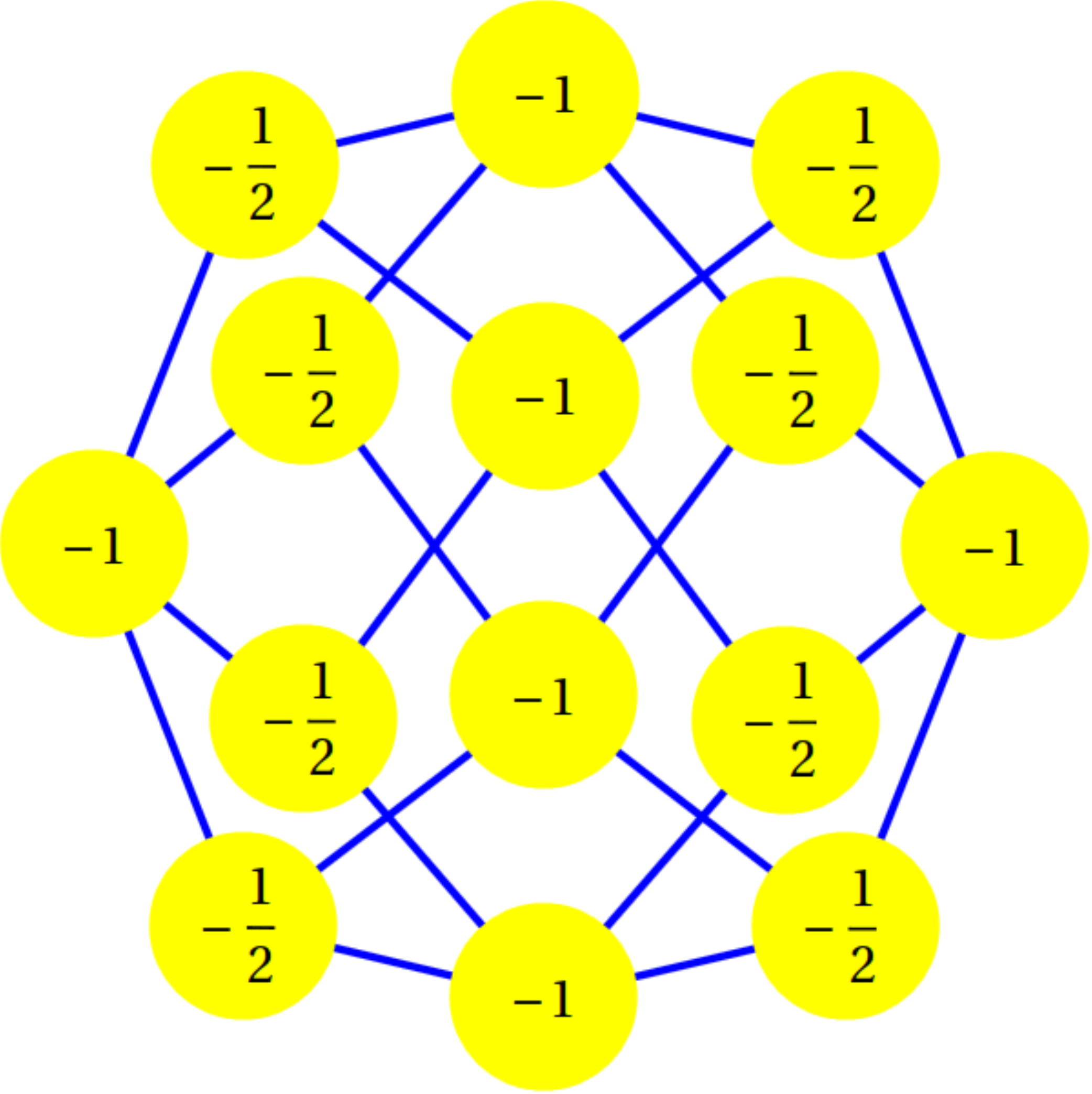}}}
\parbox{5.5cm}{\scalebox{0.25}{\includegraphics{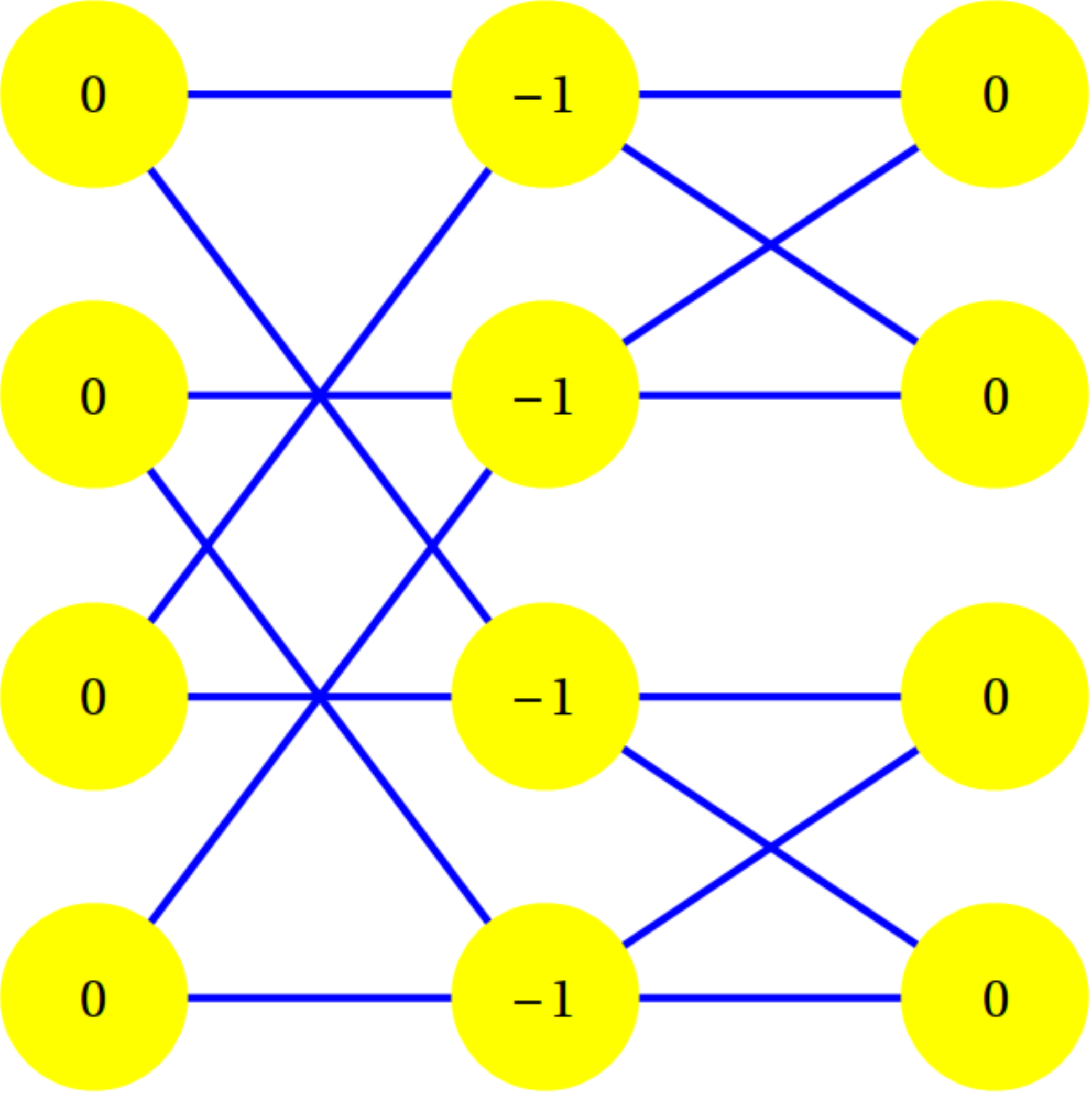}}}
}

\parbox{16.8cm}{ \parbox{5.5cm}{ Archimedean } \parbox{5.5cm}{ Catalan } \parbox{5.5cm}{ Butterfly } }

\pagebreak

\parbox{16.8cm}{
\parbox{6cm}{\scalebox{0.28}{\includegraphics{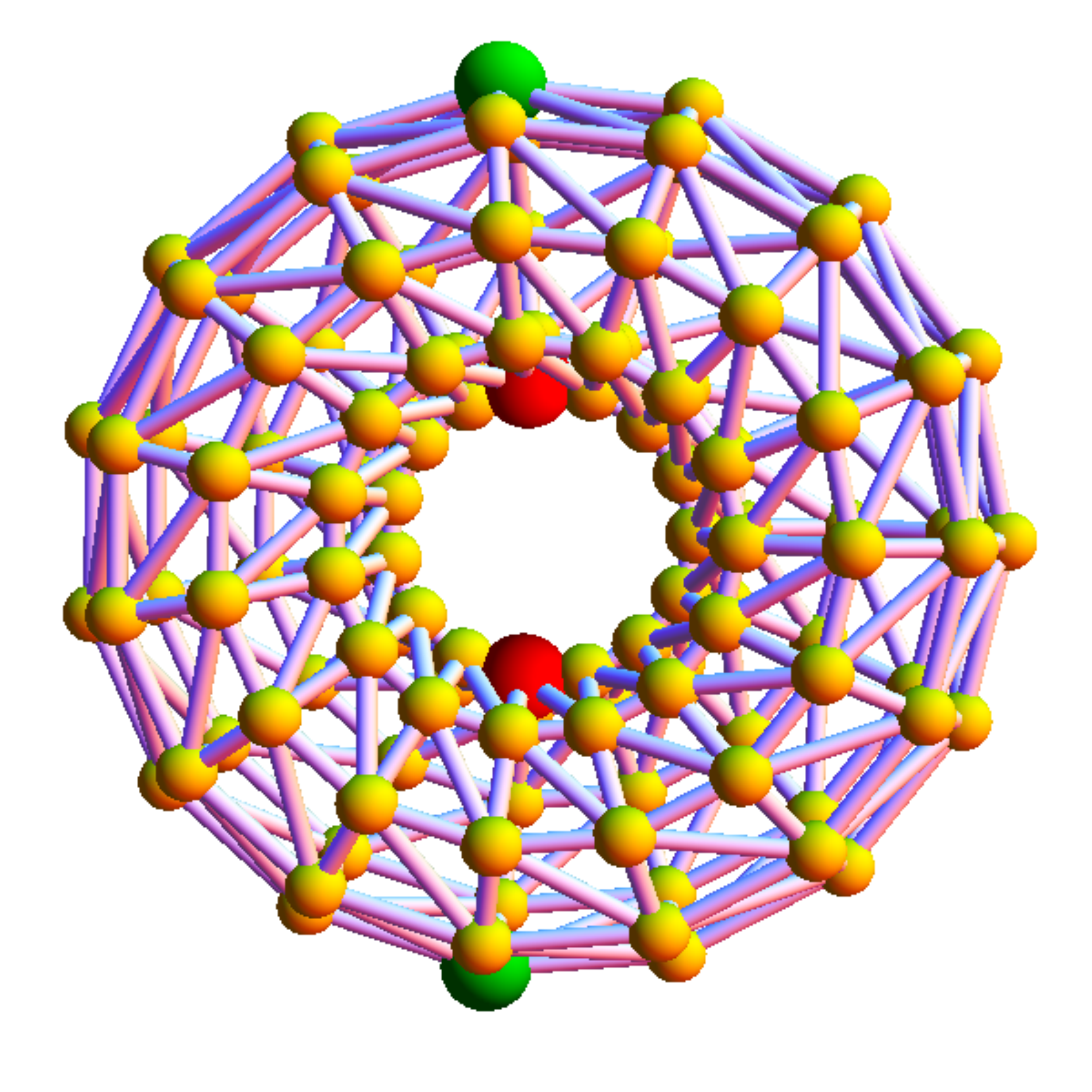}} }
\parbox{10cm}{
{\bf An example for Poincar\'e-Hopf.}
The figure shows a graph $G=(V,E)$ of Euler characteristic $0$. To get a Morse function on $V$, it is 
embedded in $3$-dimensional space by a map $r: V \to R^3$ such that all vertices have
different heights. The function $f(v) = (0,0,1) \cdot r(v)$ is then injective.
There are $4$ critical points. The minimum has index $1$ because the unit sphere is empty
of Euler characteristic $0$. The maximum $v$ has index $1$ because the unit sphere is a 
circle of Euler characteristic $0$. Therefore $\chi(S^-(v))=1-\chi(\emptyset)=1$.
The saddle points have index $-1$ because each unit sphere consist of two linear graphs
of total Euler characteristic $2$ so that $1-\chi(S^-(v))=2$ in those cases and the sum of the indices is $0$. 
}
}

\vspace{1cm}

\parbox{16.8cm}{
\parbox{6cm}{\scalebox{0.28}{\includegraphics{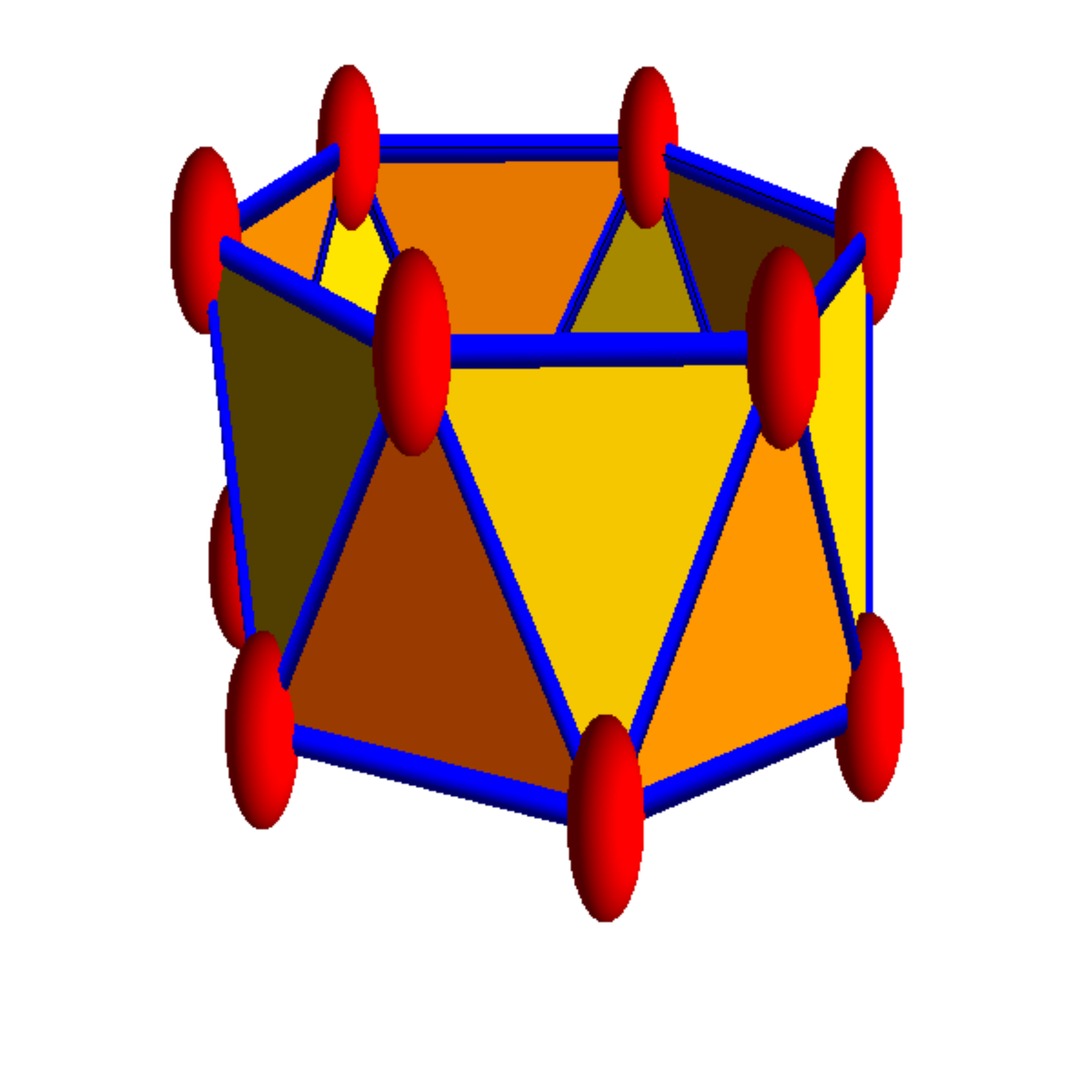}} }
\parbox{10cm}{
An orientable $2$-dimensional graph $G$ with boundary. The maximal simplices are triangles.
$G$ is a triangularization of a {\bf cylinder}. The Euler characteristic is $0$. 
It is orientable. One way to describe an orientation is to fix a permutation
on each of the triangles. The orientations are compatible, if the induced permutations on 
the edges cancel on the intersection of adjacent triangles. The sum of $df$ over the graph leads to 
cancellations, where triangles intersect. Only the boundary part survives. This leads to the
sum of $f$ along the boundary. The discrete Stokes theorem is what most mathematicians and especially
physicists have in mind, when thinking about this theorem. Writing down the theorems in the 
continuum involves only fighting with notation and language as well as battling the concept of limit but
there are no new ideas involved. Calculus flavors like nonstandard analysis
\cite{Nelson77} absorb those difficulties with additional language. While cultural and pedagogical barriers
makes this forbidding to be used for many decades to come, it is today already 
a powerful intuitive tool for many mathematicians. 
}
}

\vspace{1cm}

\parbox{16.8cm}{
\parbox{6cm}{\scalebox{0.28}{\includegraphics{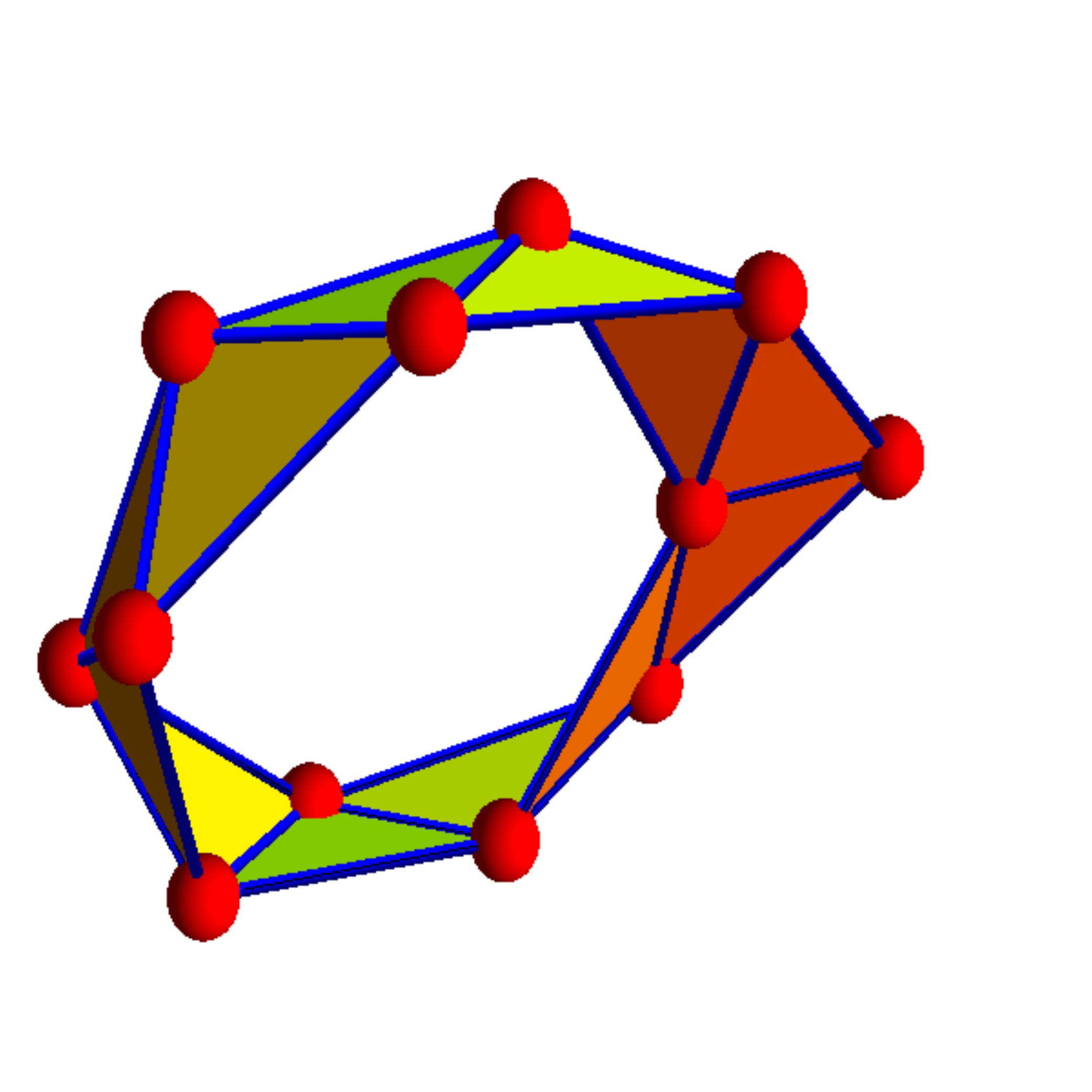}} }
\parbox{10cm}{
A non-orientable two-dimensional graph with boundary.  It is a
discrete {\bf M\"obius strip}. Also here, the Euler characteristic is zero, as 
the classical M\"obius strip has. The maximal simplices are triangles.
We can not find a $2$-form $f(x,y,z)$ which assigns the value
$1$ to each triangle and is compatible at the intersections. 
As in the continuum, we would have to look
at a double cover of the surface. In graph theory, this can be realized
by taking two copies of this M\"obius strip and then glue them together.
When realized in space this would produce a double twisted band which topologically
is a cylinder with boundary. Drawings of topologists like Alexandroff \cite{Alexandroff} 
or Fomenko \cite{Fomenko} suggest to think in terms of sceleton graphs and forget about the Euclidean fillings. 
While $n$-dimensional simplices are traditionally defined as convex sets in $R^{n+1}$, in graph theory, 
they are just complete graphs $K_{n+1}$. Graph theory is pedagogically easier to grasp because one can 
draw $K_{n+1}$ on paper, while thinking in $n$-dimensional space is a psychological 
hurdle at first. Intuition about higher dimensional space is a constant theme in 
Poincar\'e's writings (i.e. \cite{Poincare1895}). Today, we worry about this less
philosophically but pedagogically. 
}
}

\vspace{1cm}

\parbox{16.8cm}{
\parbox{6cm}{\scalebox{0.28}{\includegraphics{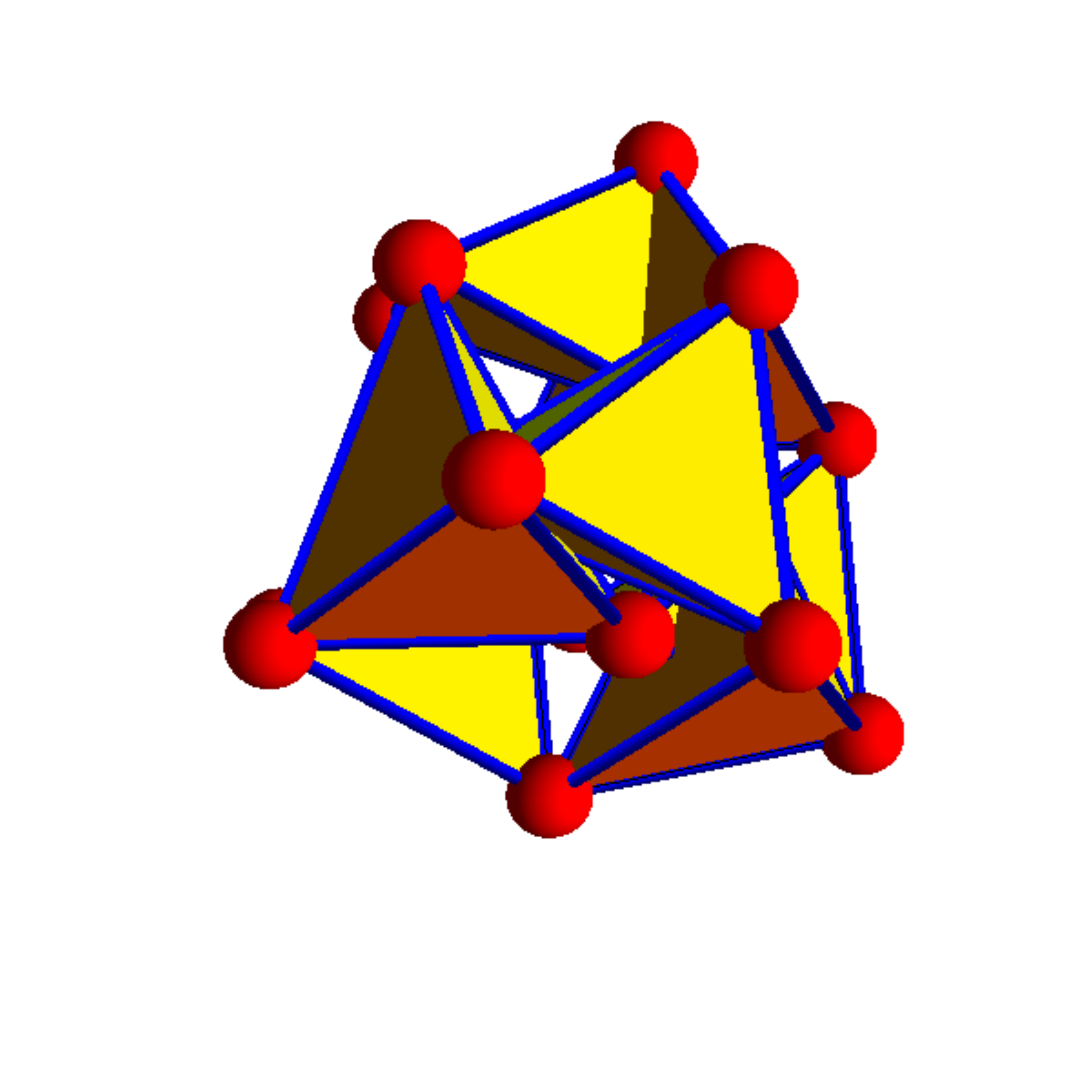}} }
\parbox{10cm}{
{\bf A three dimensional graph}.
This particular three dimensional graph $G$ has Euler characteristic $\chi(G) = -8$. 
The unit spheres $S(v)$ at each point consists of 2 triangles at each point
so that this is not a geometric graph. It is a three dimensional version of a 
Hensel type graph. Such examples of simplicial complexes appeared when struggling
to define what a "polyhedron" or "polytope" is. Because of the many incompatible
definitions inducing well known fallacies, disgusted topologists have abandoned
the term "polyhedron" and started on a clean slate with simplicial complexes. 
Today, the false proofs and refutations are well understood and the term is clear.
It is interesting however that the notion of "polyhedron" can be defined entirely
{\bf graph theoretically} without any reference to Euclidean space or an embedding in a 
topological space. A polyhedron is a graph with or without boundary which 
after a possible truncation or stellation = kising becomes a two dimensional 
graph with boundary. A tetrahedron needs a truncation
to become two dimensional. A dodecahedron needs a stellation of the 12 faces to become a 
two dimensional graph.  The graph shown to the left is not a polytope. The unit spheres are
not connected
}
}

\vspace{1cm}

\parbox{16.8cm}{
\parbox{6cm}{\scalebox{0.28}{\includegraphics{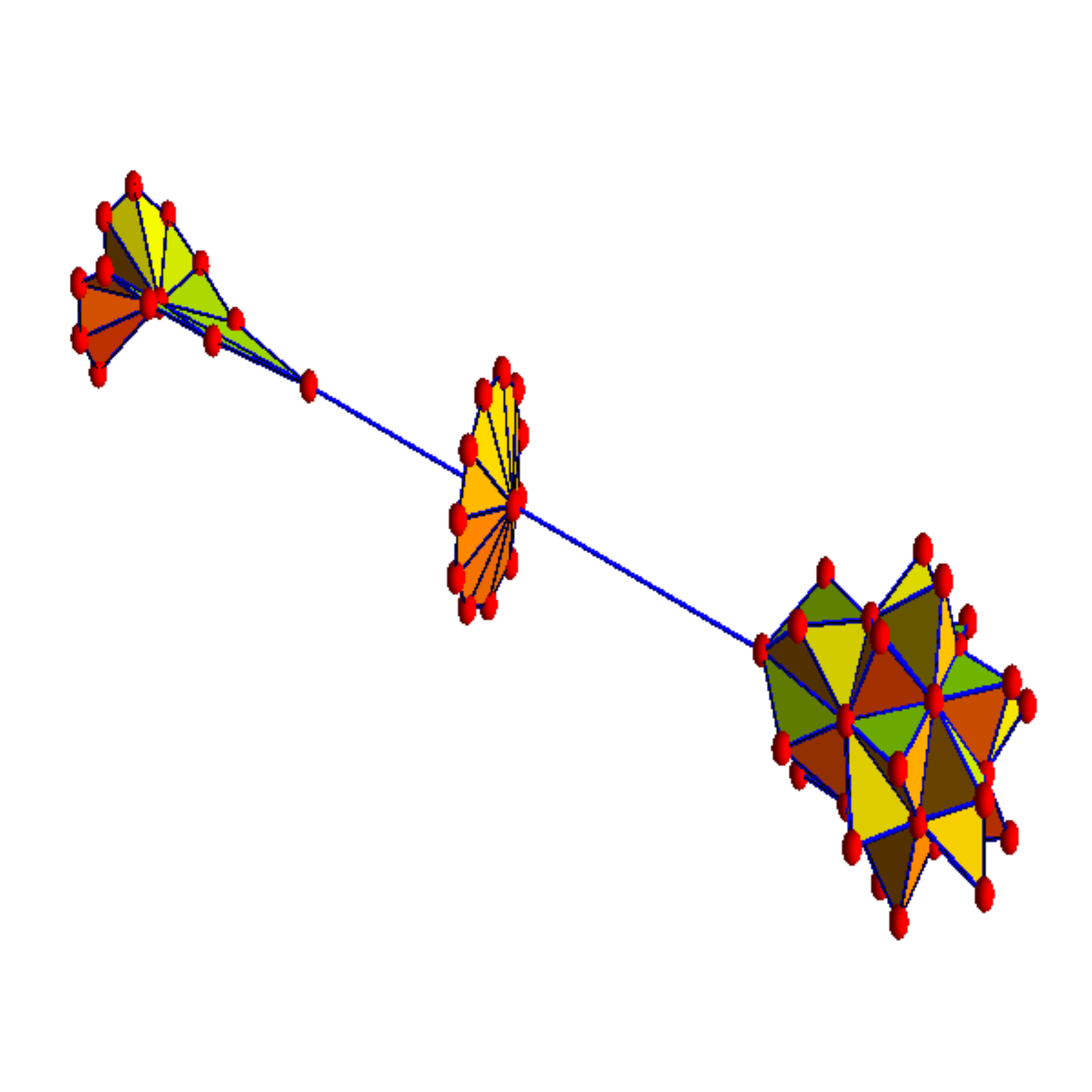}} }
\parbox{10cm}{
This graph $G$ is not a $2$-dimensional geometric graph with boundaries because some unit spheres 
have fractional dimension. The graph itself has fractional dimension.
Gauss-Bonnet and Poincar\'e-Hopf hold on this graph. For
Stokes theorem, we would have to see graphs as chains, a larger set with group structure and
which is closed under the boundary operation.
The boundary $dG$ is no more a graph but a chain, an integer-valued
function on ${\mathcal{G}} = \bigcup_k G_k$. The group elements in this group are usually
written as a chain $G=\sum_{v \in {\mathcal{G}}} a_v v$. 
Only if $dG=\sum_v b_v v$ and all $b_v$ are equal to $1$, we can forget about 
the group structure and have all information about the chain or simplicial complex 
stored as a graph. The group structure and a Euclidean realization 
can be rebuilt from this if $G$ is a graph. But for $dG$, this is not possible in general. }}

\vspace{1cm}

\parbox{16.8cm}{
\parbox{6cm}{\scalebox{0.28}{\includegraphics{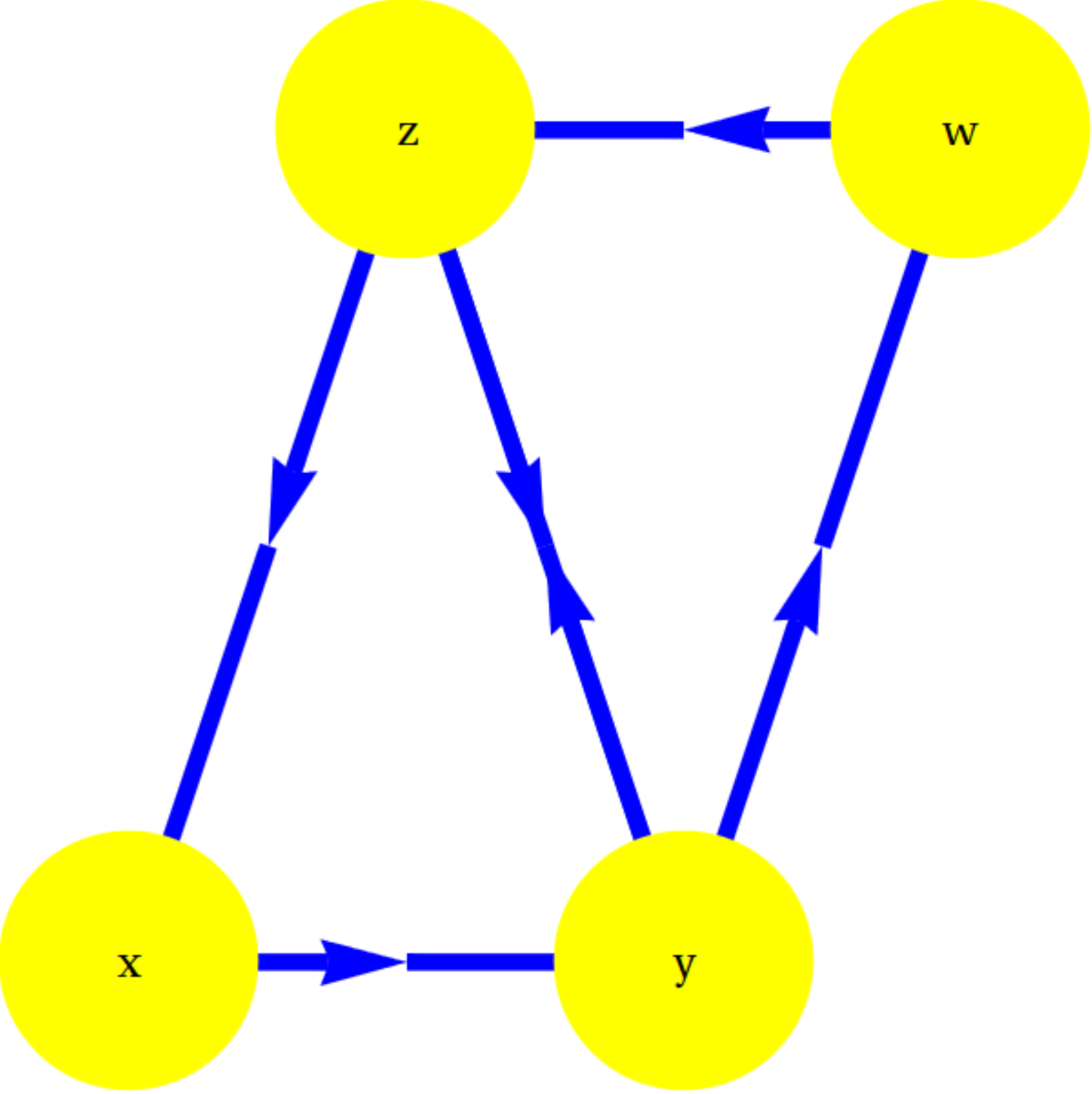}} }
\parbox{10cm}{
{\bf A simple case for Green-Stokes.}
A $2$-dimensional graph in which each point is a boundary point. 
This graph is orientable. We can assign a $2$-form $m$ which gives the value $1$ 
to the two triangles $(x,y,z),(y,w,z)$. 
The orientation defines an orientation on the boundary $d(x,y,z) = (y,z) - (x,z) + (x,y)$
and $d(y,w,z) = (w,z) - (y,z) + (y,w)$. Given a $1$-form $f$, we have
$df(x,y,z) = f(y,z) - f(x,z) + f(x,y)$ and $df(y,w,z) = f(w,z) - f(y,z) + f(y,z)$. 
We have $\sum_{G_2} df = f(y,z) - f(x,z) + f(x,y)  +  f(w,z) - f(y,z) + f(y,z) = 
- f(x,z) + f(x,y)  +  f(w,z) + f(y,z) = f(z,x) + f(x,y) +  f(w,z) + f(y,z) = \int_{dG_{1}} f$.
If $f$ is a gradient $f=dg$, then $f(x,y) = g(y)-g(x)$ and $\int_{dG} f = 
f(z,x) + f(x,y) +  f(w,z) + f(y,z) = g(x)-g(z) + g(y)-g(x) + g(z)-g(w) + g(z)-g(y)=0$. 
}
}

\vspace{1cm}

\parbox{16.8cm}{
\parbox{6cm}{\scalebox{0.28}{\includegraphics{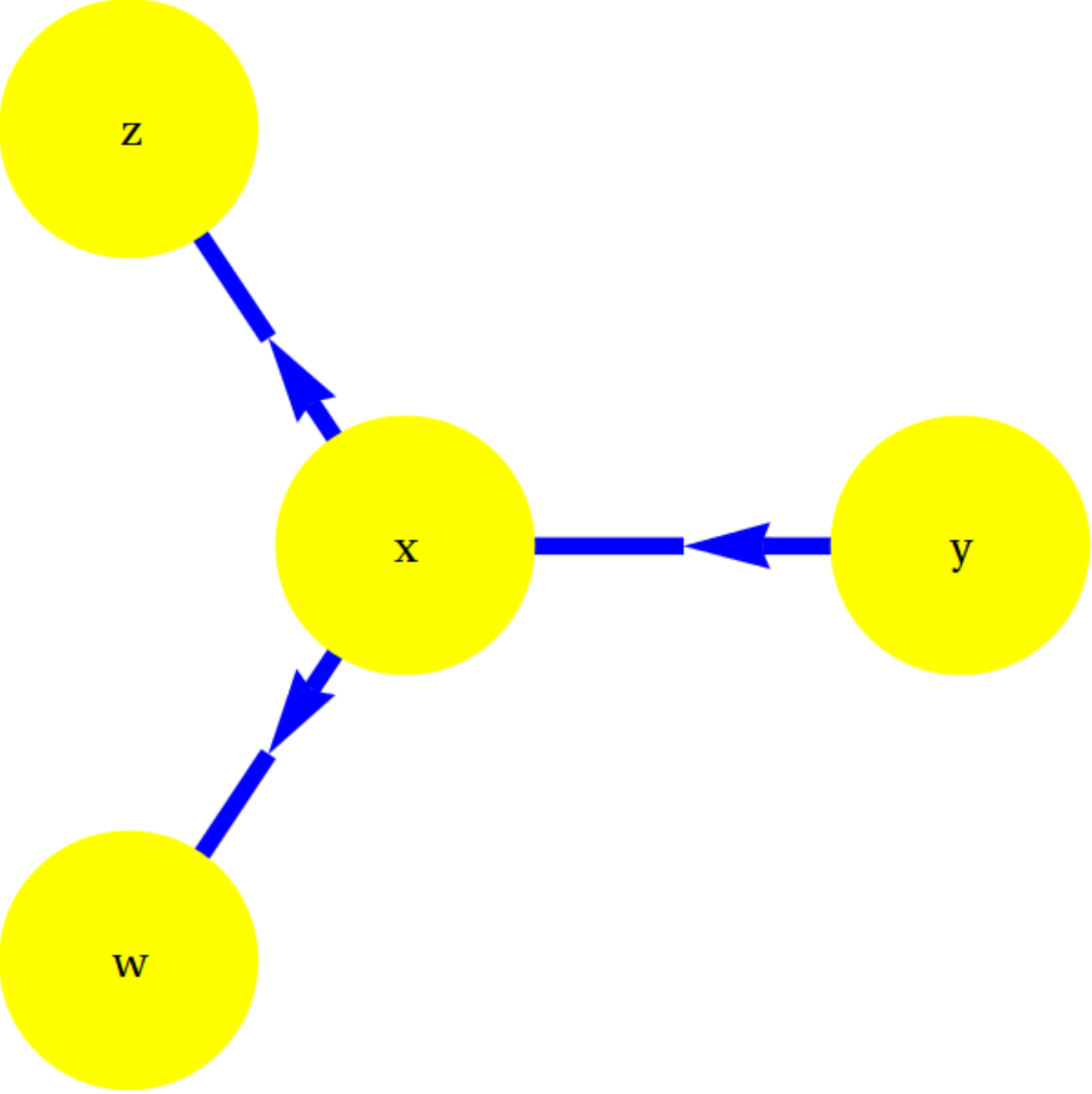}} }
\parbox{10cm}{
{\bf The boundary of a graph is in general a chain, not a graph.}
The star $S_3$ is not a "one dimensional 
graph with boundary" because $S(x) = P_3$ is not sphere like. 
The maximal simplices are the three peripheral vertices. We can attach an arbitrary orientation
but it is not compatible at the center. Given a $1$-form $f$, we have $df(x,y)=f(y)-f(x)$
and $df(x,w)=f(x)-f(w)$ and $df(z,x) = f(x)-f(z)$. The sum $\sum_{G} df=f(y)-f(w) + f(z)-f(x)$
is not the sum of $f$ over the boundary. 
The boundary of a graph is in general no more a graph. This prompted Poincar\'e 
to introduce notions like {\bf chains}, which is 
an integer valued functions on $\bigcup_k G_k$. The group of such functions on $G_k$ is called a 
{\bf simplicial $k$-chain}. Elements of the group $C_k$ are called {\bf chains} and 
can be written as $\sum_v a_v v$.  The sequence of groups $C_k \to C_{k-1}$ form a 
{\bf simplicial chain complex} and defines 
{\bf homology groups} $H_k$ whose rank $b_k$ are the {\bf homological Betti numbers}.
}
}

\vspace{1cm}

\parbox{16.8cm}{
\parbox{6cm}{\scalebox{0.28}{\includegraphics{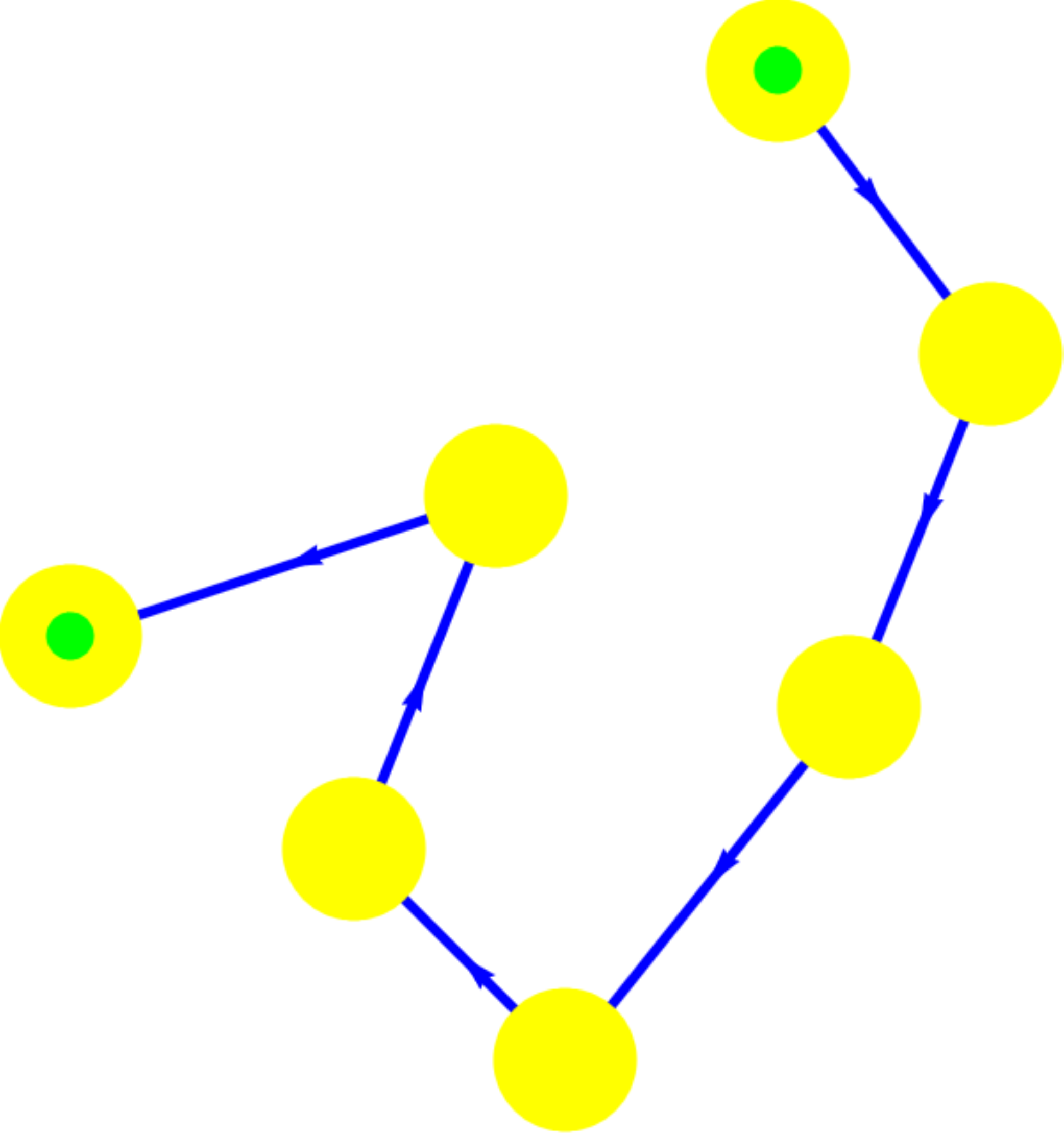}} }
\parbox{10cm}{
{\bf The fundamental theorem of line integrals.}  
A polygon without self intersections defines a one dimensional graph. If it is not
closed, it is a graph with boundary. Every connected component has two boundary points. 
The Euler characteristic of a connected component is $0$ for a closed polygon and $1$
for a polygon with boundary. The curvature is $1/2$ at each end and $0$ in the interior.
Given an injective $f$, the index $i_f$ is $1$ at local minima and $-1$ at local maxima except
at the boundary. The symmetric index $j_f$ agrees with curvature. 
The maximal simplices are edges. An {\bf orientation} on a one dimensional graph $G$ with boundary
makes it a {\bf directed graph}. The later can be imposed on
any graph but already the star graph has shown that it does not lead to an orientation even if the
graph is one dimensional. Orientation needs compatibility at intersections of maximal simplices.
Given a $0$-form $f(x)$ on vertices, the line integral $\sum_{k=0}^{n-1} df(k,k+1) = \sum_{k=0}^{n-1} f(k+1)-f(k)$ 
reduces to $f(n)-f(0) = \sum_{v \in dG} f(v)$. 
}
}

\vspace{1cm}

\parbox{16.8cm}{
\parbox{6cm}{\scalebox{0.28}{\includegraphics{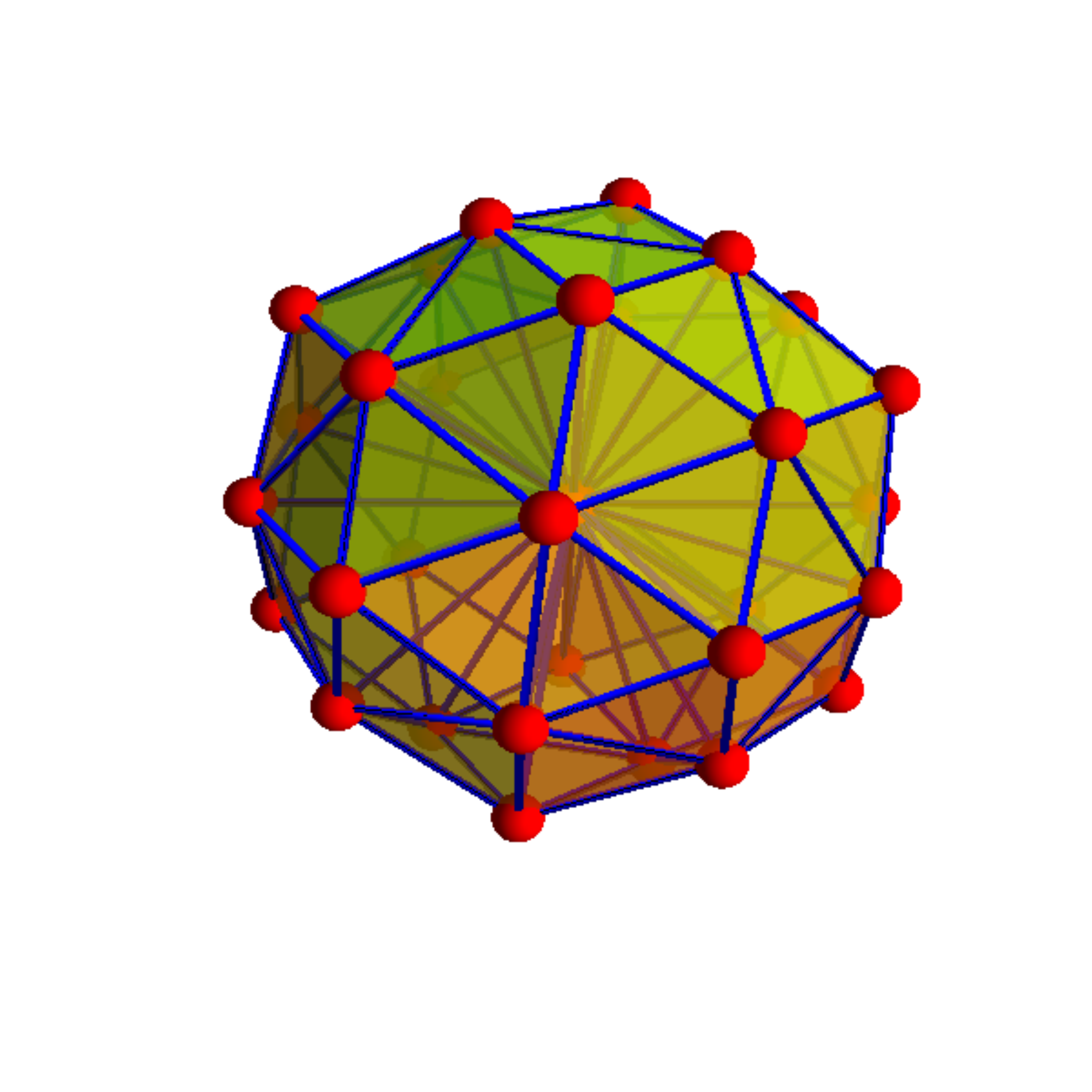}} }
\parbox{10cm}{
{\bf An example for the divergence theorem.}
A $3$-dimensional graph with boundary. It could be seen as the unit ball in a larger
$3$ dimensional graph. The boundary is the unit sphere of a vertex in that graph.
The maximal simplices are three dimensional  tetrahedra. There is an orientation on them. 
Given a $2$ form $f$ we get a $3$-form $df$. It is called
{\bf divergence} and given by $df(x,y,z,w) = f(y,z,w) - f(x,z,w) + f(x,y,w) - f(x,y,z)$.
When evaluated on all the tetrahedra, only the part on the boundary survives 
because contributions over matching faces cancel. This leads to 
the sum $\sum_{v \in dG_2} f(v)$, which is the flux of $f$ through the boundary. 
The later sum is the discrete analogue of the flux of the field $f$ through the 
boundary surface. 
}
}

\vspace{1cm}

\parbox{16.8cm}{
\parbox{6cm}{\scalebox{0.28}{\includegraphics{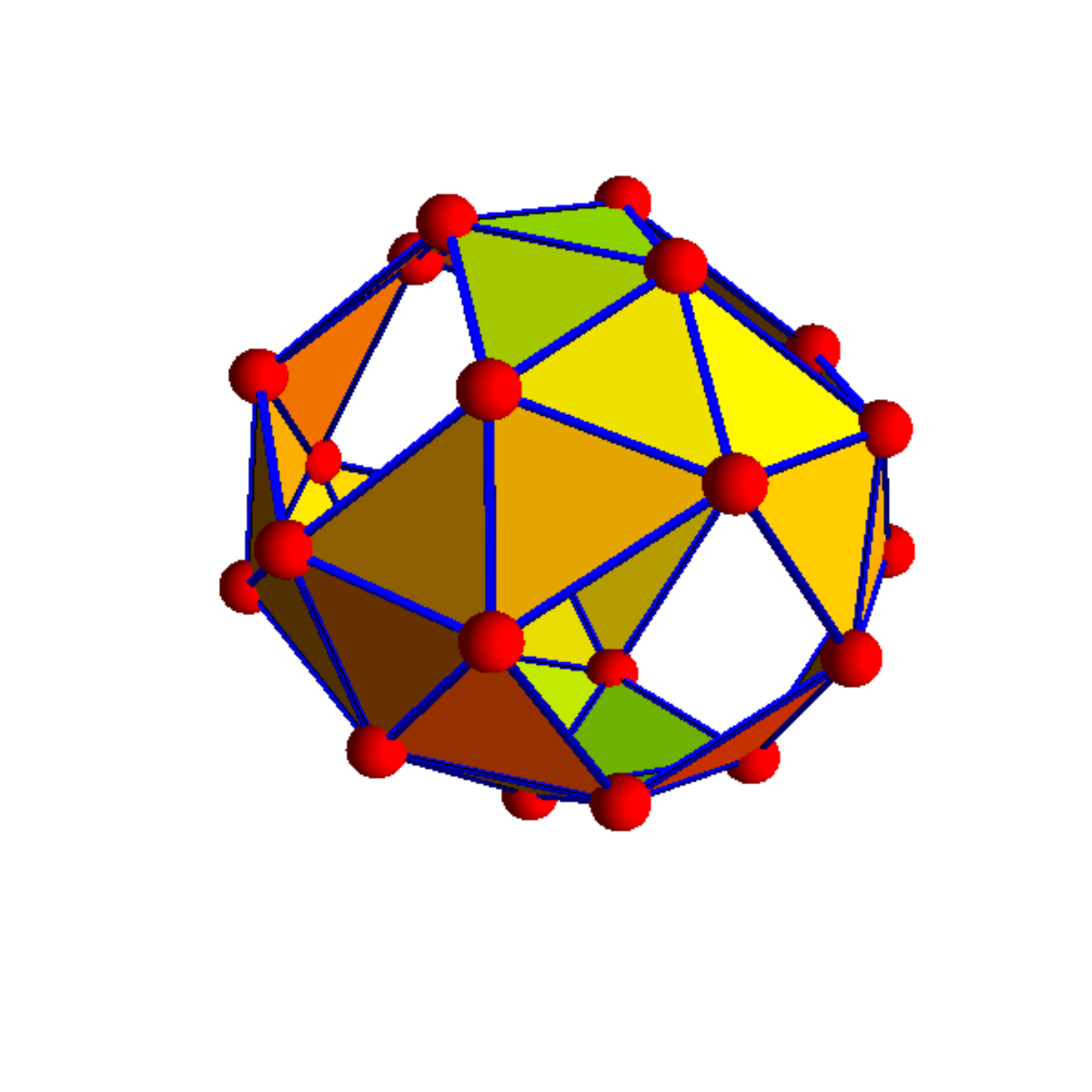}} }
\parbox{10cm}{
A case for Stokes Theorem. This two dimensional graph with boundary is a 
{\bf snub cube}. It has $|V|=24$ vertices and $|E|=60$ edges. 
The boundary consists of $6$ circular graphs $C_4$. The Euler characteristic is $-4$. 
If all $32$ triangles are oriented counterclockwise when looking from outside,
then each of the $6$ holes is oriented clockwise. Given a one form $f$ on the edges, we can 
compute the curl $df$ on each triangle and add it up. What survives is the sum of $f$
along the $6$ boundary curves. 
}
}

\vspace{0cm}

\parbox{16.8cm}{
\parbox{6cm}{\scalebox{0.28}{\includegraphics{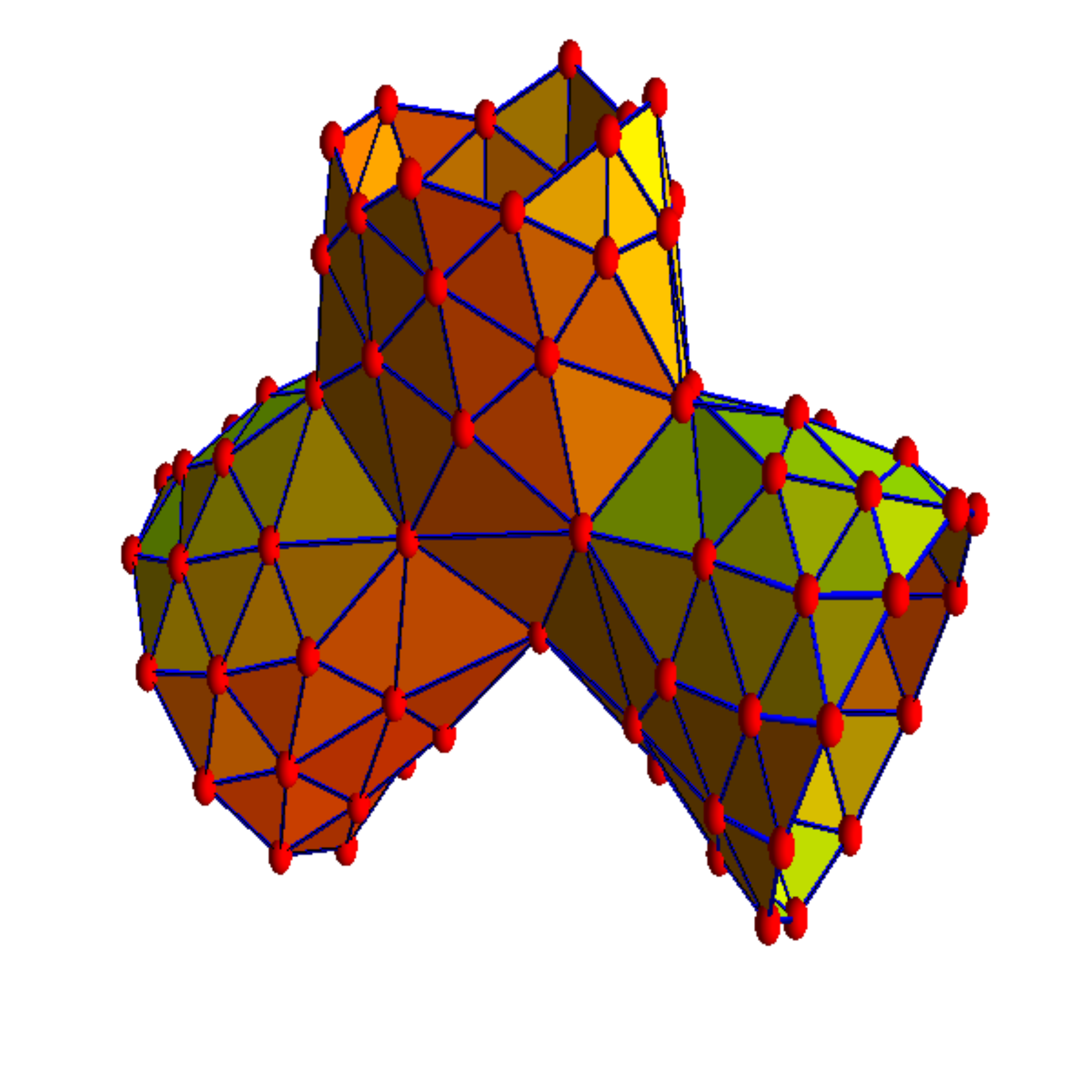}} }
\parbox{10cm}{
Poincar\'e introduced the notion of {\bf cobordism} for manifolds. 
The analogue for graphs is straightforward: 
two $k$-dimensional graphs are {\bf cobordant} if there is a 
$(k+1)$-dimensional graph $G$ with boundary, such that the boundary $dG$
has one part isomorphic to $G_1$ and an other isomorphic to $G_2$. 
This leads to the equivalence relation of {\bf h-cobordism} for graphs
which is weaker than graph isomorphism: two $k$-dimensional graphs $G_1,G_2$ are {\bf h-cobordant}
if they are cobordant with a $k+1$ dimensional graph $G=(V,E)$ and such that there is an injective
$f: V \to R$ for which at every interior point of $G$, the function $f$ restricted 
to the sphere $S(v)$ has exactly $2$ critical 
points. Any two cyclic graphs $C_n,C_m$ for example are h-cobordant for $n,m \geq 4$.
The figure illustrates a cobordism between a cyclic graph and two cyclic graphs. 
It is not a h-cobordism.}}

\vspace{1cm}

\parbox{16.8cm}{
\parbox{6cm}{\scalebox{0.28}{\includegraphics{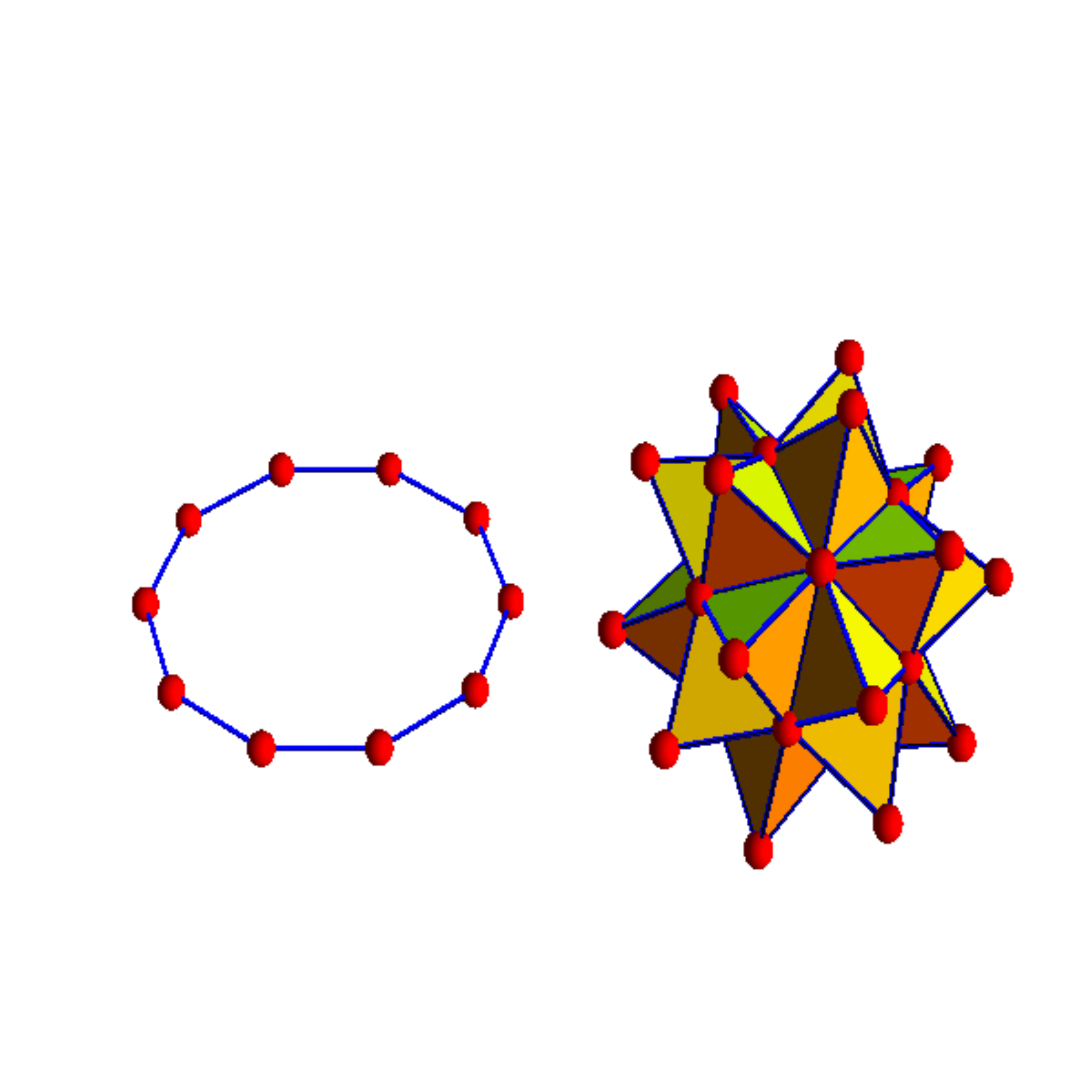}} }
\parbox{10cm}{
{\bf Reeb's theorem} motivates the following inductive definition: 
a $k$-dimensional graph $G=(V,E)$ is {\bf sphere-like } if it is an
empty graph or if every unit sphere $S(v)$ is sphere-like and the minimal
number of critical points among all injective functions $f$ on $V$ is equal to $2$ 
and inductively, $f$ restricted to each unit sphere $S(w)$ shares the same minimization 
property.  Examples: $0$-dimensional spheres are $P_2$. One-dimensional
spheres are $C_n$ for $n \geq 4$. Two-dimensional spheres
can be realized as convex polyhedra with triangular faces. In general,
a $k$-dimensional sphere like graph is a triangularization of a $k$-dimensional 
classical sphere. A function $f$ with $2$ critical points defines a filtration of
sets. Start at a minimum $v$, define  $D_0=\{v\}$ and form the strict inclusions
$D_k \subset D_{k+1} = D_k \cup \bigcup_{w \in d D_k} S^+(w)$ until $D_m = G$. 
The boundaries of $D_k,k=1,\dots, m-1$ are spheres which are by definition all 
h-cobordant for $0 \leq k<m$. The discrete analogue of Reeb's theorem is that
a {\bf sphere like graph} is a {\bf simplicial sphere}, a sceleton of 
a simplicial complex which when stuffed with Euclidean meat is homeomorphic to a sphere.}
}

\vspace{1cm}

\parbox{16.8cm}{
\parbox{6cm}{\scalebox{0.28}{\includegraphics{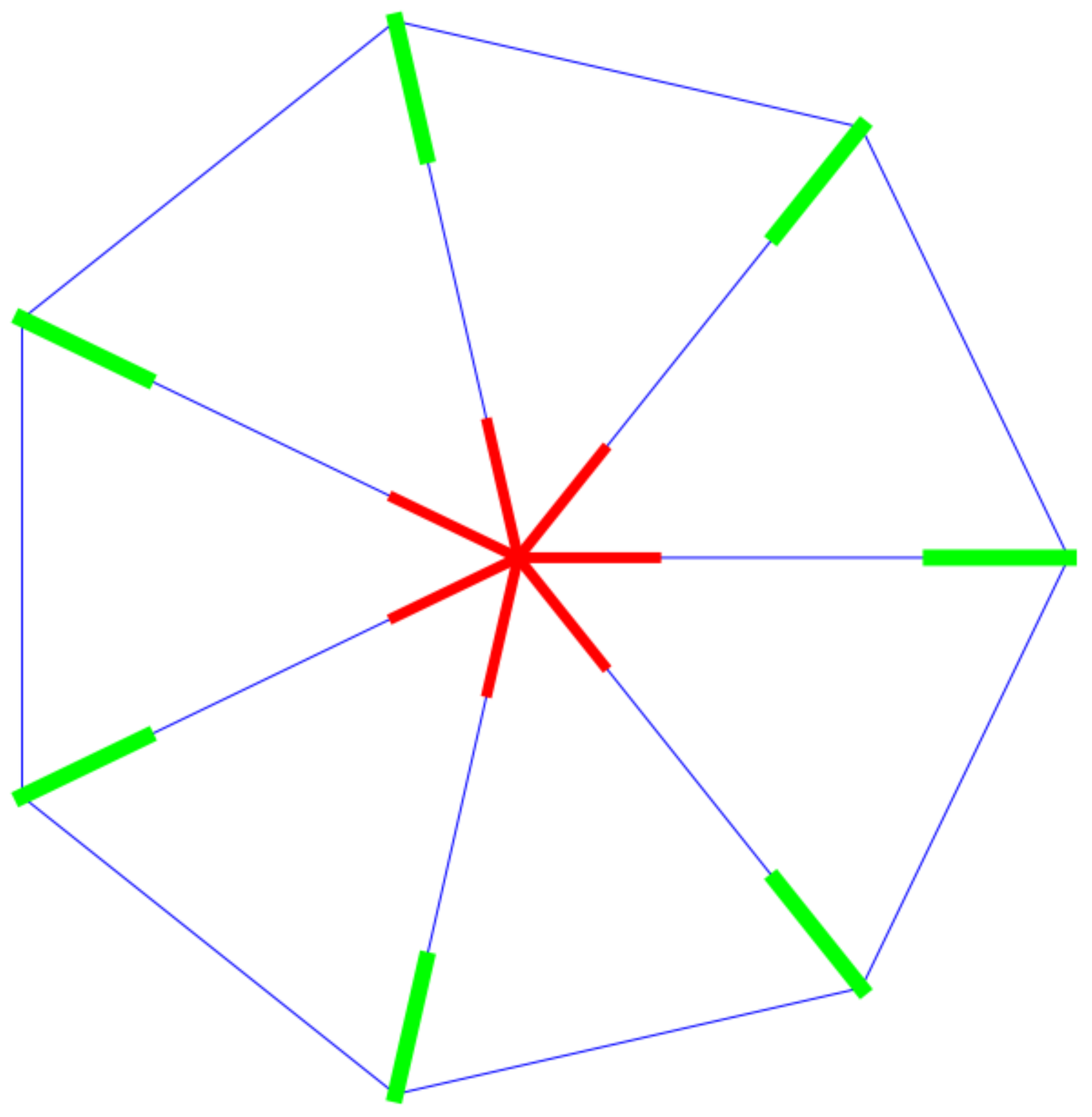}} }
\parbox{10cm}{
{\bf Proof of the transfer equations.}
The figure illustrates the proof of the transfer equations 
\fcolorbox{brightgreen}{brightgreen}{$\sum_{v \in V} V_{k-1}(v) = (k+1) v_k$}
in the case $k=1$, where the sum of the degrees over all vertices is twice the number of edges. 
If the {\bf $k$-degree} of a vertex $v$ is defined as the number of $k$-simplices
directly adjacent to $v$, then the sum of all $k$-degrees is $k+1$ times the number
of $k$ simplices. Draw connections from every vertex to 
every center of any $k$-simplex. We can count these connections in two 
different ways: by first summing up all leading a given vertex and then sum
over $V$ or by seeing that every simplex has $k+1$ connections and then sum 
over the simplices. 
}
}

\vspace{1cm}

\parbox{16.8cm}{
\parbox{6cm}{\scalebox{0.28}{\includegraphics{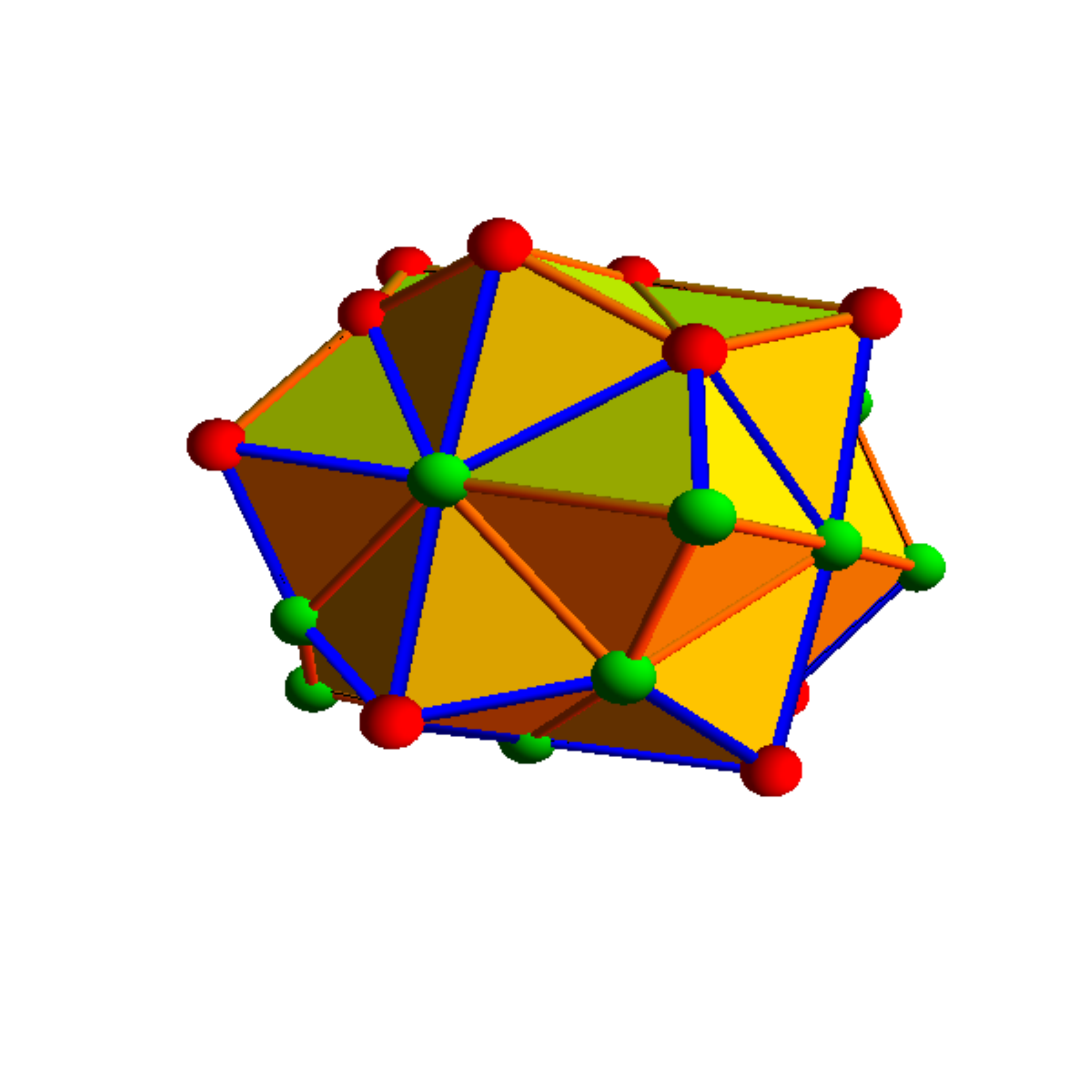}} }
\parbox{10cm}{
{\bf The $W_k$ in the intermediate equations.} An injective function $f$ on $V$
defines the sets $W_k(v)$ of all $k$ simplices in the sphere $S(v)$
for which at least one vertex $w$ satisfies $f(w)<f(v)$ and an other vertex satisfies $f(u)>f(v)$.
The figure shows the unit sphere of a point $v$ in a three dimensional graph. The vertices $w$ of the
sphere are colored differently according to whether they are in $S^-(v)$ or $S^+(v)$. 
The set $W_1(v)$ consists of edges in $S(v)$ which connect vertices with different colors. 
The intermediate equations tell that 
$\sum_{v \in V} W_1(v) = v_{2}$ is the number of the triangles. Indeed, each edge in $W_1(v)$ 
defines together with $v$ a triangle. This map from $\bigcup_{v \in V} W_1(v)$ to $G_2$ 
is a bijection as triangles can not appear twice in the listing: there is exactly one vertex
in a triangle where $f(w)$ is between the values of the two others. 
}
}

\vspace{1cm}

\parbox{16.8cm}{
\parbox{6cm}{\scalebox{0.28}{\includegraphics{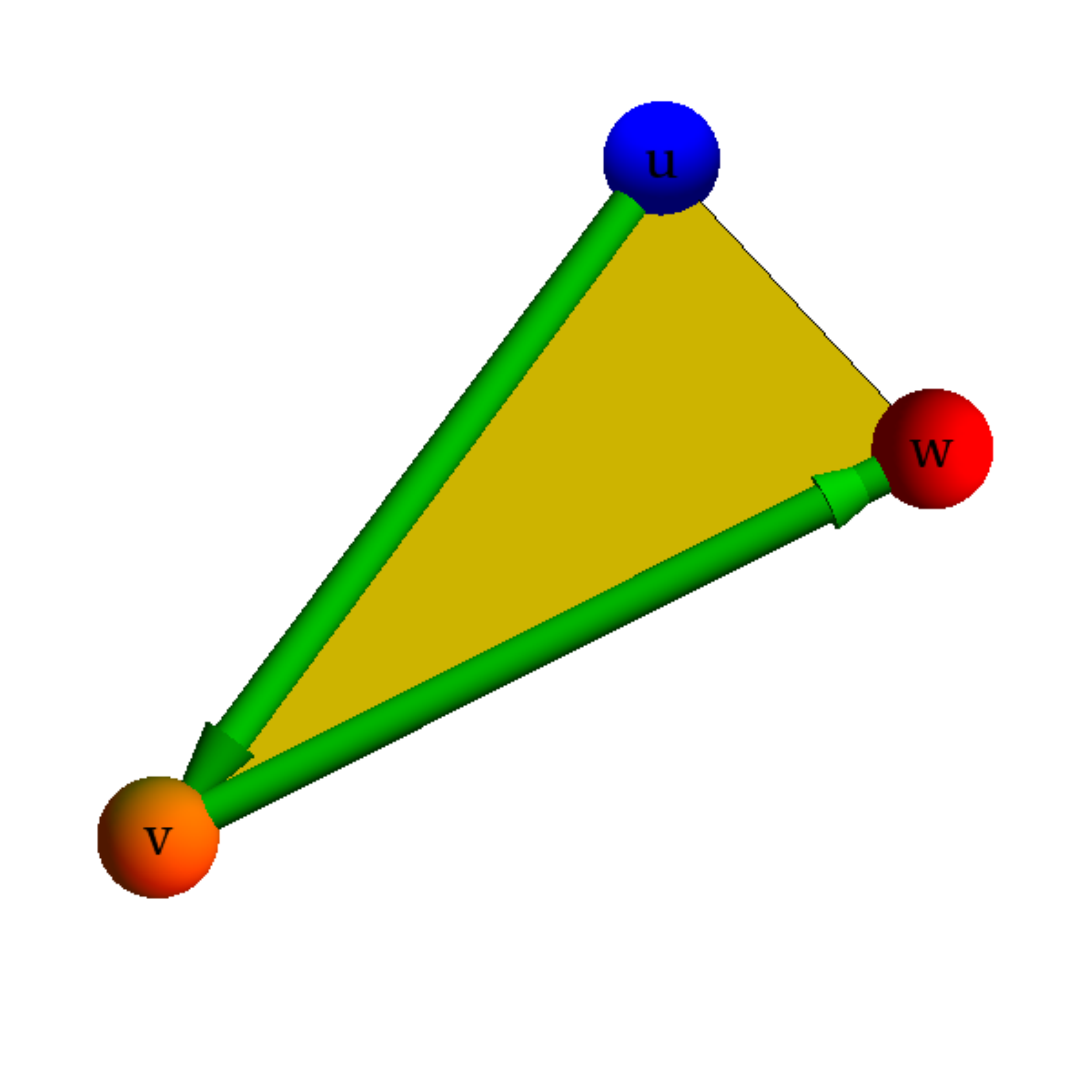}} }
\parbox{10cm}{
{\bf Proof of the intermediate equations.} Given an injective function $f$ on $V$
and a vertex $v$, we can look at the set $W_k(v)$ of all $k$ simplices in the sphere $S(v)$
for which at least one vertex $w$ satisfies $f(w)<f(v)$ and an other vertex satisfies $f(u)>f(v)$. 
The figure illustrates the intermediate equations
\fcolorbox{brightgreen}{brightgreen}{$\sum_{v \in V} W_k(v) = k v_{k+1}$}
in the case $k=1$. If we look the ordering given by the values of $f$ and start at $w$, 
we get to $v$ and then to $u$. 
The ordering introduced through $f$ makes this triangle unique in the sense
that the triangle is not counted from any of the other $2$ vertices.
In the case $k=2$, we look at the $v_{k+1}$ tetrahedra. In each tetrahedron, there are $2$ vertices
which neighbor vertices with larger and smaller values. Each of these vertices $v$ has a triangle in the
unit sphere which adds to the $W_k(v)$ sum. This leads to $\sum_{v \in V} W_2(v) = 2 v_{3}$.
}
}

\pagebreak 

\section{Historical links}

For the history of the classical Stokes theorem, see \cite{Katz79}. 
The history of topology \cite{Dieudonne1989,HistoryTopology}.
The collection \cite{HistoryTopology} contains in particular an article on the history of graph theory. 
A story about Euler characteristic and polyhedra is told in \cite{Richeson}. 
For an introduction to Gauss-Bonnet with historical pointers to early discrete approaches see \cite{Wilson}.
For Poincar\'e-Hopf, the first volume of \cite{Spivak1999} or \cite{Hirsch}. For Morse theory and 
Reeb's theorem \cite{Mil65}.
Poincar\'e proved the index theorem in chapter VIII of \cite{poincare85}. Hopf extended it to 
arbitrary dimensions in \cite{hopf26}.
Gauss-Bonnet in higher dimensions was proven first independently by Allendoerfer \cite{Allendoerfer} 
and Fenchel \cite{Fenchel} for surfaces in Euclidean space and extended jointly by 
Allendoerfer and Weil \cite{AllendoerferWeil} to closed Riemannian manifolds. 
Chern gave the first intrinsic proof in \cite{Chern44}.

\vfill

\end{document}